# Rigorous validation of a Hopf bifurcation in the Kuramoto-Sivashinsky PDE

Jan Bouwe van den Berg, Elena Queirolo

September 30, 2020

**Abstract**

We use computer-assisted proof techniques to prove that a branch of non-trivial equilibrium solutions in the Kuramoto–Sivashinsky partial differential equation undergoes a Hopf bifurcation. Furthermore, we obtain an essentially constructive proof of the family of time-periodic solutions near the Hopf bifurcation. To this end, near the Hopf point we rewrite the time periodic problem for the Kuramoto–Sivashinsky equation in a desingularized formulation. We then apply a parametrized Newton-Kantorovich approach to validate a solution branch of time-periodic orbits. By contruction, this solution branch includes the Hopf bifurcation point.

## 1 Introduction

The goal of this paper is to give a proof of a Hopf bifurcation found numerically in the Kuramoto-Sivashinsky partial differential equation (PDE) in one spatial dimension with $2\pi$-periodic boundary conditions in space:

$$\begin{cases} \partial_t u = -\partial_x^4 u - \gamma \partial_x^2 u + \partial_x(u^2), \\ u(t,x) = u(t, x+2\pi). \end{cases} \quad (1)$$

Here $\gamma > 0$ acts as the bifurcation parameter. The appealing qualities of the PDE (1) are its relative simplicity and the presence of a low order nonlinear term. Despite its simplicity, it exhibits complex dynamics, and it serves as a model for weak turbulence in laminar flows and more generally as a paradigm for studying spatiotemporal chaos. We refer to [6], [7] and [8] for more background material on the Kuramoto-Sivashinsky equation.

In this paper we prove the existence of a Hopf bifurcation from a nonhomogeneous stationary solution to (1) by *a posteriori* validation of numerical computations. The Kuramoto–Sivashinsky equation has been well studied from the perspective of validated numerics. Early results include [16], where stationary solutions of (1) are found and validated. In [12] rigorous numerics is used to prove the existence of complicated, chaotic trajectories in the ordinary differential equation (ODE) describing steady states of the Kuramoto–Sivashinsky equation. Much more recently, by a powerful combination of a rigorous integrator for parabolic PDEs and topological arguments, the existence of an invariant set with chaotic dynamics was proven [13] for the semiflow induced by the PDE (1). Furthermore, the bifurcation diagram of stationary solutions of the Kuramoto-Sivashinsky equation and their stability was studied in [1], see also [15] for the pitchfork bifurcation problem. We refer to [2], [14], [4] and [5] for validation of periodic orbits in Kuramoto–Sivashinsky. Here we complement these results by an analysis of the Hopf bifurcation problem. Indeed, the periodic solutions proven in [2, 4, 5, 14] bifurcate from an inhomogenous stationary state, but a mathematically rigorous analysis of this Hopf bifurcation problem has not previously been undertaken.

To tackle this, we generalize our previous desingularization work on (Hopf) bifurcations in ODEs [9] to the PDE setting. The main strength of this approach is that validating the existence of a Hopf bifurcation and the branch of time-periodic solutions emanating from it, is reformulated as a regular branch following problem. The crucial advantage is that the solution branch of periodic orbits does not become singular at the bifurcation point, hence continuation can be carried out using rigorous computer-assisted pseudo-arclength continuation techniques. There is then a straightforward correspondence between solution of the desingularized system on the one hand and those of the original system on the other. This "blow-up" technique complements the approaches to computer-assisted bifurcation analysis introduced in [1] and [15]. While we focus on Hopf bifurcations in this paper, our technique is well-suited to study more general symmetry breaking bifurcations from nonhomogeneous stationary states.

To outline the desingularization approach, it is convenient to change variables. We consider a time $L$-periodic ($L$ being unknown for the moment), space $2\pi$-periodic solution $u(t,x)$ to (1). We apply a time rescaling $\tau = \frac{2\pi t}{L}$ such that $u(\tau, x)$ is time $2\pi$-periodic. The rescaled equation is written as

$$\begin{cases} \partial_\tau u = -\lambda_1 \partial_x^4 u - \lambda_1 \lambda_2 \partial_x^2 u + \lambda_1 \partial_x(u^2), \\ u(\tau, x) = u(\tau, x + 2\pi), \end{cases} \qquad (2)$$

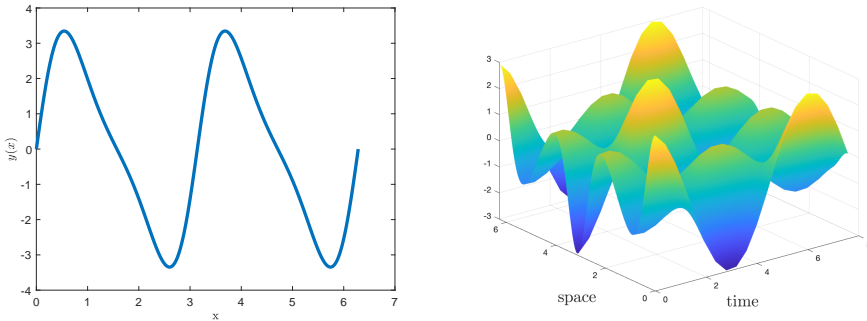

Figure 1: On the left: validated $2\pi$-periodic stationary solution $y(x)$ near the Hopf bifurcation point . On the right: rescaled time and space periodic profile $z(\tau, x)$ near the Hopf bifurcation.

with $\lambda_1 = \frac{L}{2\pi}$ and $\lambda_2 = \gamma$.

If a time-dependent solution $u(\tau, x)$ is close enough to the Hopf bifurcation, it is natural to rewrite $u$ as the sum of a stationary solution $y(x)$ to (2) and a time-dependent perturbation $az(\tau, x)$, where $z$ is of order 1 and $a \in \mathbb{R}$ is the (small) amplitude of the perturbation. The precise formulation can be found in Section 2.3. Following this rewriting, we validate a branch of solutions $(\lambda_1, \lambda_2, a, y, z)$. If the amplitude $a$ changes sign along this branch, then at $a = 0$ a Hopf bifurcation occurs: the periodic solutions for $a \neq 0$ bifurcate from the stationary state $y$ at $a = 0$.

The validated solution close to the bifurcation is presented in Figure 1, where the stationary component $y(x)$ of the solution is presented on the left and the perturbation profile $z(\tau, x)$ generated by the Hopf bifurcation is plotted on the right. The full solution of (1) is the sum of the stationary component and the perturbation multiplied by the amplitude $a$, which tends to zero at the bifurcation point. Thanks the the validated numerics presented in this paper, we can prove that the Kuramoto-Sivashinsky equation (1) undergoes a Hopf bifurcation at the parameter value $\lambda_2^* = \gamma^* = 0.298776358114 + [-r, r]$ with validation radius $r = 8.63 \cdot 10^{-11}$.

By adding a continuation condition, we can prove that $a' \neq 0$, where the derivative is taken with respect a "parametrization" variable of the solution branch. This guarantees the "non-degenacy" of the Hopf bifurcation in the sense that we prove a branch of non-constant periodic orbits which converge to the equilibrium state. In this way we avoid the need for a center manifold analysis. Furthermore, for this result we do not need to check that the conjugate eigenvalue pair, which we find at $a = 0$ to lie on the imaginary axis, crosses from one side of the imaginary axis to the other side (when following the branch of stationary solution in $\gamma$). Such an analysis requires an additional computer-assisted proof, which would not be too difficult, but we did not pursue it in this paper. Finally, we remark that in order to establish sub- or supercriticality of the bifurcation, one could determine the sign of $a''$. In the ODE setting, this analysis has been performed in great generality in [9], and we are confident that such an approach would work in the PDE setting as well. Since working out the details is somewhat tedious and the current paper is already quite heavy on technical detail, we decided not to include it here.

The details of the result presented graphically in Figure 1 are discussed in Section 9. The main challenge in transforming the ODE methods in [9] to the PDE setting of the current paper are two-fold. First, the two-dimensional Fourier series lead to diagonally dominant operators, whose diagonal depends on the *a priori* unknown time period $\lambda_1$. This requires additional estimates to perform the tail (truncation) estimates. Furthermore, the nonlinearity includes a derivative, which can be controlled because it is lower order than the bi-Laplacian, but again this requires additional bounds. Both of these are discussed in Section 6. We stress that lifting the techniques from ODEs to PDEs is *not* difficult from a technical point of view. One may find this disappointing, but we view it as the pivotal strength of the approach. Indeed, the blow up approach presented is extremely general in its simplicity. It seems probable that the bounds developed in this paper can be applied with minor modifications to many other parabolic PDEs with periodic (or Neuman) boundary conditions.

The outline of the paper is as follows. In Section 2 we introduce the functional analytic setup, which is analogous to that of [5] and [4], employing Fourier series and a convenient choice of the norm in Fourier space. The quantitative contraction theorem of Newton-Kantorovich type, which underlies the computer-assisted proof, is presented in Section 3. Sections 4 and 5 are devoted to introducing and bounding the relevant projections and other linear operators. The necessary results for analyzing the tails of the Fourier modes are collected in Section 6. In Section 7 the "Newton-Kantorovich" bounds for single periodic orbits are derived. Subsequently, these are generalized to bounds for the continuation case in Section 8. Finally, in Section 9 we collect all ingredients to prove the existence result discussed above.

## 2 The Kuramoto-Sivashinsky equation and the Hopf bifurcation system

### 2.1 The Kuramoto-Sivashinsky equation

In Section 1, we rescaled the Kuramoto-Sivashinsky equation to

$$\partial_t v = -\lambda_1 \lambda_2 \partial_x^4 v - \lambda_1 \partial_x^2 v + 2\lambda_1 \partial_x(v^2), \tag{3}$$

so that we can study solutions which are $2\pi$-periodic in space and time.

To write $v$ in terms of its Fourier expansion in $t$ and $x$ we introduce the two-dimensional bi-infinite sequence $u = (u_{jk})_{j,k\in\mathbb{Z}}$, with $u_{jk} \in \mathbb{C}$, such that

$$v(t,x) = \sum_{j,k\in\mathbb{Z}} u_{jk} e^{\mathrm{i}(jt+kx)}. \tag{4}$$

From here onwards, we will consistently use $j$ for time related indices and $k$ for space related indices. We introduce the notation

$$\mathcal{F}v = u$$

for the two-dimensional Fourier transform, with $u$ defined above.

We recall some standard properties of Fourier sequences that will be used extensively in the following. The convolution of two two-dimensional bi-infinite sequences $u$ and $w$ is defined by

$$(u * w)_{jk} = \sum_{j_1+j_2=j} \sum_{k_1+k_2=k} u_{j_1 k_1} w_{j_2 k_2}.$$

Furthermore it holds that

$$u * w = \mathcal{F}(\mathcal{F}^{-1}(u) \cdot \mathcal{F}^{-1}(w)),$$

where $\mathcal{F}^{-1}$ is the inverse Fourier transform and the multiplication $\cdot$ is considered pointwise in $t$ and $x$. The derivative operators satisfy

$$\mathcal{F}(\partial_t v(t,x)) = (\mathrm{i} j u_{jk})_{j,k\in\mathbb{Z}} \quad \text{and} \quad \mathcal{F}(\partial_x v(t,x)) = (\mathrm{i} k u_{jk})_{j,k\in\mathbb{Z}}.$$

We denote the multiplication of a two-dimensional bi-infinite sequence by its index in the time or space direction by the operator J or K, respectively, defined as

$$\mathrm{J} : u \mapsto (j u_{jk})_{j,k\in\mathbb{Z}} \tag{5}$$

and
$$\mathrm{K} : u \mapsto (ku_{jk})_{j,k\in\mathbb{Z}}. \tag{6}$$

We rewrite Equation (3) in the space of Fourier sequences as
$$\mathrm{i}\mathrm{J}u = -\lambda_1\lambda_2\mathrm{K}^4 u + \lambda_1 \mathrm{K}^2 u + \mathrm{i}\lambda_1 \mathrm{K}(u*u), \tag{7}$$

and we introduce the zero finding problem
$$f(\lambda_1, \lambda_2, u) \stackrel{\text{def}}{=} \mathrm{i}\mathrm{J}u + \lambda_1\lambda_2 \mathrm{K}^4 u - \lambda_1 \mathrm{K}^2 u - \mathrm{i}\lambda_1 \mathrm{K}(u*u) = 0. \tag{8}$$

**Remark 2.1.** *Writing explicitly*
$$f_{jk}(\lambda_1, \lambda_2, u) = \mathrm{i}j u_{jk} + \lambda_1\lambda_2 k^4 u_{jk} - \lambda_1 k^2 u_{jk} - \lambda_1 \mathrm{i} k (u*u)_{jk}, \qquad j,k \in \mathbb{Z},$$
*we notice that* $f_{00}(\lambda_1, \lambda_2, u) \equiv 0$.

When $u$ is a solution to (8), any shift in time or space of $u$ is also a solution. In order to lift this invariance, we add two phase conditions to our problem of the form
$$\langle \mathrm{i}\mathrm{K}\tilde{u}, u \rangle = 0 \quad \text{and} \quad \langle \mathrm{i}\mathrm{J}\tilde{u}, u \rangle = 0, \tag{9}$$

where $\tilde{u}$ is an approximation of the solution and
$$\langle u, w \rangle = \sum_{j,k\in\mathbb{Z}} u_{jk} w_{jk},$$

for $u, w \in \mathbb{C}^{\mathbb{Z}^2}$. Conditions (9) are the Fourier sequence representations of
$$\int_0^{2\pi} v \tilde{v}_x dx = 0 \quad \text{and} \quad \int_0^{2\pi} v \tilde{v}_t dt = 0,$$

with $\tilde{v}$ a numerically determined approximation of the solution. These are standard conditions to fix shifts (see e.g. [10]).

The full problem then becomes
$$\begin{cases} f_{jk}(\lambda_1, \lambda_2, u) = 0, & \text{for all } (j,k) \in \mathbb{Z}^2 \setminus \{(0,0)\}, \\ \langle \mathrm{i}\mathrm{K}\tilde{u}, u \rangle = 0, \\ \langle \mathrm{i}\mathrm{J}\tilde{u}, u \rangle = 0, \end{cases} \tag{10}$$

and one expects a one-dimensional solution branch, due to the two parameters, two scalar equations and considering Remark 2.1.

## 2.2 The functional spaces $\ell^1_{\nu_2}$ and $\mathcal{L}^1_\nu$

In this section, we introduce the spaces of one and two-dimensional bi-infinite sequences. Following [3], for one-dimensional bi-infinite sequences we introduce the norm

$$\|u\|_{\ell^1_{\nu_2}} \stackrel{\text{def}}{=} \sum_{k \in \mathbb{Z}} \nu_2^{|k|} |u_k|, \tag{11}$$

and the space of one-dimensional exponentially decreasing sequences as

$$\ell^1_{\nu_2} = \{u = (u_k)_{k \in \mathbb{Z}}, u_k \in \mathbb{C} \mid \|u\|_{\ell^1_{\nu_2}} < \infty\}.$$

The convolution of two sequences $y$ and $w$ in $\ell^1_{\nu_2}$ is defined as

$$(y * w)_k = \sum_{k_1 + k_2 = k} y_{k_1} w_{k_2}.$$

The norm and space introduced for bi-infinite sequences are very similar. For a given vector $\nu = (\nu_1, \nu_2)$, with $\nu_i \geq 1$ for $i = 1, 2$, we introduce the $\mathcal{L}^1_\nu$-norm in the space of bi-infinite sequences

$$\|u\|_{\mathcal{L}^1_\nu} \stackrel{\text{def}}{=} \sum_{j,k \in \mathbb{Z}} \nu_1^{|j|} \nu_2^{|k|} |u_{jk}|, \tag{12}$$

that is the two-dimensional version of the norm already used in [10]. It resembles closely the norm used in other papers discussing Kuramoto-Sivashinsky, such as [4] and [5]. In particular, the $\mathcal{L}^1_\nu$-norm has similar properties to the $M$-norm presented in [4]. We define the space of two-dimensional exponentially decreasing sequences as

$$\mathcal{L}^1_\nu = \{u = (u_{jk})_{j,k \in \mathbb{Z}}, u_{j,k} \in \mathbb{C} \mid \|u\|_{\mathcal{L}^1_\nu} < \infty\}.$$

The $\mathcal{L}^1_\nu$ space shares a lot of properties with the $\ell^1_{\nu_2}$ space. First, $\mathcal{L}^1_\nu$ is a Banach space, and

$$\|u * v\|_{\mathcal{L}^1_\nu} \leq \|u\|_{\mathcal{L}^1_\nu} \|v\|_{\mathcal{L}^1_\nu}.$$

Furthermore, there is a trivial embedding of $\ell^1_{\nu_2}$ in $\mathcal{L}^1_\nu$

$$E_\ell : \ell^1_{\nu_2} \to \mathcal{L}^1_\nu$$
$$y \mapsto E_\ell y : (E_\ell y)_{jk} = \delta_{j0} y_k,$$

which satisfies $\|y\|_{\ell^1_{\nu_2}} = \|E_\ell y\|_{\mathcal{L}^1_\nu}$. Thanks to this embedding, we can make sense of the convolution between an element $y$ of $\ell^1_{\nu_2}$ and an element $u$ of $\mathcal{L}^1_\nu$, defining it as

$$(y * u)_{jk} \stackrel{\text{def}}{=} (E_\ell y * u)_{jk} = \sum_{k_1 + k_2 = k} y_{k_1} u_{jk_2}.$$

## 2.3 The Hopf bifurcation

A sequence $y \in \ell^1_{\nu_2}$ is a stationary solution of (8) if it satisfies

$$F_1(\lambda_2, y) \stackrel{\text{def}}{=} \lambda_2 \mathrm{K}^4 y - \mathrm{K}^2 y - \mathrm{i}\mathrm{K}(y * y) = 0. \tag{13}$$

We remark that $(F_1)_0 = 0$ is trivially satisfied. To avoid shifts in space of the solution to (13), as we did for the time-dependent solution in Section 2.1, we add the phase condition

$$\langle \mathrm{i}\mathrm{K}\tilde{y}, y \rangle_\ell \stackrel{\text{def}}{=} \langle \mathrm{i} E_\ell \mathrm{K}\tilde{y}, E_\ell y \rangle = \sum_{k \in \mathbb{Z}} \mathrm{i} k \tilde{y}_k y_k = 0, \tag{14}$$

with $\tilde{y}$ an approximation of the solution.

Similar to the Hopf bifurcation validation presented in [9], we use a blow up approach in the neighborhood of a Hopf bifurcation. To this end, we rewrite the solution $u \in \mathcal{L}^1_\nu$ to (10) in the neighborhood of a Hopf bifurcation as the sum of the stationary solution $y \in \ell^1_{\nu_2}$ to (13) and a perturbation $az$, with amplitude $a \in \mathbb{R}$ and $z \in \mathcal{L}^1_\nu$ of order 1, that is

$$\begin{cases} u = E_\ell y + az, \\ \langle \mathrm{J}^2 \overline{\tilde{z}}, z \rangle = 1, \end{cases} \tag{15}$$

with $\tilde{z}$ a numerically determined approximation of $z$. The second equation ensures that $z$ does not tend to a (singular) stationary orbit as the Hopf bifurcation point is approached.

Replacing $u$ in (10) by $E_\ell y + az$, we arrive at

$$\begin{cases} \langle \mathrm{J}^2 \overline{\tilde{z}}, z \rangle = 1, \\ \langle \mathrm{i}\mathrm{K}(E_\ell \tilde{y} + \tilde{a}\tilde{z}), (E_\ell y + az) \rangle = 0, \\ \langle \mathrm{i}\mathrm{J}(E_\ell \tilde{y} + \tilde{a}\tilde{z}), (E_\ell y + az) \rangle = 0, \\ \mathrm{i}\mathrm{J}(E_\ell y + az) + \lambda_1 \lambda_2 \mathrm{K}^4(E_\ell y + az) \\ \quad - \lambda_1 \mathrm{K}^2(E_\ell y + az) - \lambda_1 \mathrm{i}\mathrm{K}((E_\ell y + az) * (E_\ell y + az)) = 0. \end{cases} \tag{16}$$

We add the equations (13) and (14) and use this information on $y$ to simplify the terms involving only $y$ in the last equation. Concerning the phase condition in time, we recall that $y$ is time-independent, hence $\mathrm{J} y = 0$. Following a division by $a$, we consider the zero-finding problem

$$\begin{cases} \langle \mathrm{J}^2 \bar{\tilde{z}}, z \rangle - 1 = 0, \\ \langle \mathrm{i}\mathrm{K}\tilde{y}, y \rangle_\ell = 0, \\ \langle \mathrm{i}\mathrm{K}(E_\ell \tilde{y} + \tilde{a}\tilde{z}), z \rangle = 0, \\ \langle \mathrm{i}\mathrm{J}\tilde{z}, z \rangle = 0, \\ \lambda_2 \mathrm{K}^4 y - \mathrm{K}^2 y - \mathrm{i}\mathrm{K}(y * y) = 0, \\ \mathrm{i}\mathrm{J}z + \lambda_1 \lambda_2 \mathrm{K}^4 z - \lambda_1 \mathrm{K}^2 z - \lambda_1 \mathrm{i}\mathrm{K}(2 E_\ell y * z + a z * z) = 0. \end{cases} \qquad (17)$$

**Remark 2.2.** *The equation $\langle \mathrm{i}\mathrm{K}(E_\ell \tilde{y} + \tilde{a}\tilde{z}), z \rangle = 0$ is not the only choice for the phase condition on $z$. Other choices are possible, such as setting $\langle \mathrm{i}\mathrm{K}\tilde{z}, z \rangle = 0$. In particular, the choice made in this paper is unfit in case $y$ is a homogeneous steady state, because in that case the equation would be trivially satisfied at the bifurcation point. However, in the case at hand $y$ is space-dependent, hence the stated condition is appropriate.*

Let $\hat{x} = (\hat{\lambda}_1, \hat{\lambda}_2, \hat{a}, \hat{y}, \hat{z})$ be a numerical approximation of the solution of (17), with $\hat{y}$ and $\hat{z}$ having just a finite number of non-zero Fourier coefficients. We use the numerical approximation to fix the phase conditions. The problem can then be written as

$$\begin{cases} G^{\hat{x}}(\lambda_1, \lambda_2, a, y, z) = \begin{pmatrix} \langle \mathrm{J}^2 \bar{\hat{z}}, z \rangle - 1 \\ \langle \mathrm{i}\mathrm{K}\hat{y}, y \rangle_\ell \\ \langle \mathrm{i}\mathrm{K}(E_\ell \hat{y} + \hat{a}\hat{z}), z \rangle \\ \langle \mathrm{i}\mathrm{J}\hat{z}, z \rangle \end{pmatrix} = 0, \\ F(\lambda_1, \lambda_2, a, y, z) = \begin{pmatrix} F_1(\lambda_2, y) \\ F_2(\lambda_1, \lambda_2, a, y, z) \end{pmatrix} = 0. \end{cases} \qquad (18)$$

Here we have split the problem into a linear problem $G^{\hat{x}}$, where we emphasize the dependency on the numerical approximation, and a nonlinear problem $F$, with $F_1$ as defined in (13) and

$$F_2(\lambda_1, \lambda_2, a, y, z) \stackrel{\mathrm{def}}{=} \mathrm{i}\mathrm{J}z + \lambda_1 \lambda_2 \mathrm{K}^4 z - \lambda_1 \mathrm{K}^2 z - \lambda_1 \mathrm{i}\mathrm{K}(2 E_\ell y * z + a z * z).$$

We will at times use the notation $(F_1)_k$ and $(F_2)_{jk}$ to indicate the $k$-th and $(j,k)$-th Fourier coefficients of $F_1(x)$ and $F_2(x)$ respectively, where the argument will be clear from context.

The set of variables for (18) has the form $(\lambda_1, \lambda_2, a, y, z)$ with $\lambda_1, \lambda_2, a \in \mathbb{C}$, $y \in \ell^1_{\nu_2}$ and $z \in \mathcal{L}^1_\nu$. The full space of variables is thus

$$X \stackrel{\text{def}}{=} \mathbb{C}^3 \times \ell^1_{\nu_2} \times \mathcal{L}^1_\nu, \tag{19}$$

equipped with the product norm

$$\|x\| = \|(\lambda_1, \lambda_2, a, y, z)\| \stackrel{\text{def}}{=} \max\{|\lambda_1|, |\lambda_2|, |a|, \|y\|_{\ell^1_{\nu_2}}, \|z\|_{\mathcal{L}^1_\nu}\}.$$

At times, we will also use the so-called component norm

$$\|x\|_c \stackrel{\text{def}}{=} (\max\{|\lambda_1|, |\lambda_2|, |a|\}, \|y\|_{\ell^1_{\nu_2}}, \|z\|_{\mathcal{L}^1_\nu}) \in \mathbb{R}^3. \tag{20}$$

**Definition 2.3.** *The conjugate $x^*$ of $x \in X$ is defined by*

$$x^* \stackrel{\text{def}}{=} (\bar{\lambda}_1, \bar{\lambda}_2, \bar{a}, y^*, z^*),$$

*with the conjugate of $y$ and $z$ defined by*

$$\begin{aligned} y^*_k &\stackrel{\text{def}}{=} \bar{y}_{-k}, \quad k \in \mathbb{Z}, \\ z^*_{jk} &\stackrel{\text{def}}{=} \bar{z}_{-j-k}, \quad j, k \in \mathbb{Z}. \end{aligned}$$

*Let $Y$ be a space on which conjugate symmetry is defined, such as $X$, $\ell^1_{\nu_2}$, or $\mathcal{L}^1_\nu$. Then the subspace of conjugate symmetric elements is denoted by*

$$\mathcal{S}(Y) \stackrel{\text{def}}{=} \{x \in Y : x = x^*\}.$$

**Definition 2.4.** *A map $T : X \to X$ is said to preserve conjugate symmetry if $T(\mathcal{S}(X)) \subseteq \mathcal{S}(X)$.*

For $\nu' < \nu$, we have

$$\begin{pmatrix} G^{\hat{x}} \\ F \end{pmatrix} : \mathbb{C}^3 \times \ell^1_{\nu_2} \times \mathcal{L}^1_\nu \to \mathbb{C}^4 \times \ell^1_{\nu'_2} \times \mathcal{L}^1_{\nu'}$$
$$(\lambda_1, \lambda_2, a, y, z) \mapsto (v', y', z').$$

Furthermore, $(F_1)_k$ and $(F_2)_{jk}$ vanish for $k = 0$ and $(j, k) = (0, 0)$, that is $y'_0 = 0$ and $z'_{00} = 0$. Therefore we can compose with the map

$$\tilde{\Pi} : \mathbb{C}^n \times \ell^1_{\nu'_2} \times \mathcal{L}^1_{\nu'} \to \mathbb{C}^{n-2} \times \ell^1_{\nu'_2} \times \mathcal{L}^1_{\nu'}$$
$$(v', y', z') \mapsto (v'', y'', z'')$$

with $v''_i = v'_i$ for $i = 1, \ldots, n-2$ and

$$y''_k = \begin{cases} v'_{n-1} & \text{if } k = 0, \\ y'_k & \text{if } k \neq 0, \end{cases} \quad \text{and} \quad z''_{jk} = \begin{cases} v'_n & \text{if } (j,k) = (0,0), \\ z'_{jk} & \text{if } (j,k) \neq (0,0). \end{cases}$$

In this way we remove the trivial identities from the formulation and arrive at

$$H^{\hat{x}}(x) \stackrel{\text{def}}{=} \tilde{\Pi}\begin{pmatrix} G^{\hat{x}} \\ F \end{pmatrix}(x) = 0, \qquad (21)$$

which generically has a one-dimensional solution set. We abuse notation by not specifying $n = 4$ in the definition of $\tilde{\Pi}$, so that we can reuse the notation later in a slightly more general setting.

Let $x = (\lambda_1, \lambda_2, a, y, z) \in X$ and $w = (\beta_1, \beta_2, b, r, s) \in X$, then we extend the notation $\langle \cdot, \cdot \rangle$ to $X$ by defining

$$\langle x, w \rangle \stackrel{\text{def}}{=} \lambda_1 \beta_1 + \lambda_2 \beta_2 + ab + \langle y, r \rangle_\ell + \langle z, s \rangle.$$

In the following, we will assume we have computed two numerical solutions $\hat{x}_0 = (\hat{\lambda}_{10}, \hat{\lambda}_{20}, \hat{a}_0, \hat{y}_0, \hat{z}_0) \in \mathcal{S}(X)$ and $\hat{x}_1 = (\hat{\lambda}_{11}, \hat{\lambda}_{21}, \hat{a}_1, \hat{y}_1, \hat{z}_1) \in \mathcal{S}(X)$ of (18). Assuming the two numerical solutions are close to each other, i.e. $\|\hat{x}_0 - \hat{x}_1\| \ll 1$, it is reasonable to expect that the segment connecting them is an approximation of a nearby solution branch. To construct a one-to-one correspondence between the point on the segment $[\hat{x}_0, \hat{x}_1]$ and the solution branch $x(s)$, $s \in [0, 1]$, we follow [5] and add to (18) the so-called continuation equation

$$\boldsymbol{E}_s(x) = \langle x, q_s \rangle - c_s, \qquad (22)$$

where

$$q_s = q_0 + s(q_1 - q_0), \quad \text{and} \quad c_s = c_0 + s(c_1 - c_0),$$

with

$$q_{\mathsf{s}} = (v_{\mathsf{s}} + v_{\mathsf{s}}^*)/2, \quad \text{and} \quad c_{\mathsf{s}} = \langle \hat{x}_{\mathsf{s}}, q_{\mathsf{s}} \rangle, \quad \text{for } \mathsf{s} = 0, 1, \qquad (23)$$

where $v_{\mathsf{s}}$ is approximately in $\ker(DH^{\hat{x}_{\mathsf{s}}}(\hat{x}_{\mathsf{s}}))$, with $H^{\hat{x}_{\mathsf{s}}}(x)$ defined in (21). Furthermore, we assume $q_{\mathsf{s}}$, for $\mathsf{s} = 0, 1$, has finitely many non-zero Fourier coefficients, as do $\hat{x}_0$ and $\hat{x}_1$. The continuation problem is then

$$H_s(\lambda_1, \lambda_2, a, y, z) = \tilde{\Pi}\begin{pmatrix} \boldsymbol{E}_s(\lambda_1, \lambda_2, a, y, z) \\ G^{\hat{x}_s}(\lambda_1, \lambda_2, a, y, z) \\ F(\lambda_1, \lambda_2, a, y, z) \end{pmatrix}. \qquad (24)$$

**Remark 2.5.** *It follows from the definition of $H_s$ in (24) that*

$$H_s(x^*) = H_s(x)^*,$$

*where the codomain of $H_s$ is $X' = \mathbb{C}^3 \times \ell^1_{\nu'_2} \times \mathcal{L}^1_{\nu'}$ for $\nu' < \nu$. Therefore, $H_s$ preserves conjugate symmetry.*

We are interested in solving the continuation problem

$$H_s : X \to X', \qquad H_s(x) = 0 \quad \text{for all } s \in [0, 1]. \tag{25}$$

**Remark 2.6.** *A solution branch $x(s)$ of (25) goes through a Hopf bifurcation if there exist $s^*, s_0$ and $s_1$ such that $s^* \in (s_0, s_1)$, $a(s^*) = 0$ and $a(s) \neq 0$ for $s \in [s_0, s_1]\setminus\{s^*\}$. The latter statement guarantees that the time-dependent term $az$ does not vanish unless $a = 0$.*

*By continuity of $a(s)$, if it holds that $a(s_0) \cdot a(s_1) < 0$ then there exists a $s^* \in (s_0, s_1)$ such that $a(s^*) = 0$. We are thus left with proving that $a(s)$ changes sign along the validated solution branch and $a'(s^*) \neq 0$, in order to prove the existence of a Hopf bifurcation of the Kuramoto-Sivashinsky equation.*

**Remark 2.7.** *Although we have formulated the problem in the context of general pseudo-arclength continuation, we can choose the continuation equation (22) to conveniently simplify the verification in Remark 2.6. Indeed, the choice*

$$q_\mathbf{s} = (0, 0, 1, 0, 0) \in X \qquad \text{and} \qquad c_\mathbf{s} = \hat{a}_\mathbf{s}$$

*in essence turns the general (pseudo-arclength) continuation into parameter continuation with respect to $a$. In this case, if a solution exists, we immediately conclude that $a'(s) \neq 0$, provided $\hat{a}_0 \neq \hat{a}_1$. As a consequence, if $\hat{a}_0 \cdot \hat{a}_1 < 0$, we can prove that $a(s)$ changes sign along the validated curve of solutions and Remark 2.6 implies the existence of a Hopf bifurcation.*

## 3 The radii polynomial approach

For the parameter dependent zero finding problem (25), we aim to validate the segment between two numerical solutions $\hat{x}_0$ and $\hat{x}_1$, which are such that

$$H_\mathbf{s}(\hat{x}_\mathbf{s}) \approx 0, \quad \hat{x}_\mathbf{s} \in \mathcal{S}(X) \quad \text{for } \mathbf{s} = 0, 1.$$

For $s \in [0, 1]$, we introduce the notation

$$\hat{x}_s = \hat{x}_0 + s(\hat{x}_1 - \hat{x}_0).$$

Following [5], we transform the zero finding problem (25) into a parameter dependent fixed point problem

$$T_s(x) = x,$$

with $s \in [0,1]$ and $T_s$ the Newton-Kantorovich operator defined by

$$T_s : X \to X \tag{26}$$
$$x \mapsto x - A_s H_s(x), \tag{27}$$

where $A_s$ is an approximation of $DH_s(\hat{x}_s)^{-1}$ for all $s$ in $[0,1]$. The linear operator $A_s$ will be constructed in Section 6 satisfying the following conditions:

$$\begin{aligned} & A_s H_s(x) \in X \text{ for } x \in X \\ & A_s \text{ is injective for all } s \in [0,1]. \end{aligned} \tag{28}$$

In order to prove that $T_s$ is a contraction in a neighborhood of $\hat{x}_s$, we first introduce

$$B_r(x) \stackrel{\text{def}}{=} \{y \in X : \|y - x\|_X \leqslant r\}.$$

To prove that $T_s$ is a contraction in $B_r(\hat{x}_s)$ for some $r$, we determine bounds

$$Y_s \geqslant \max_{s \in [0,1]} \|T_s(\hat{x}_s) - \hat{x}_s\|_X, \tag{29}$$

and

$$Z_s(r) \geqslant \max_{s \in [0,1]} \sup_{b,c \in B_1(0) \subset X} \|DT_s(\hat{x}_s + rb)rc\|_X. \tag{30}$$

With $Y_s$ and $Z_s(r)$ computed rigorously, we can then apply the following parametrized Newton-Kantorovich theorem and validate our numerical approximations.

**Theorem 3.1.** *Let $A_s$ satisfy the condition (28), let $T_s : X \to X$ be defined as in (26), let $Y_s$ and $Z_s(r)$ be defined as in (29) and (30) respectively. Let the radii polynomial be defined as*

$$p(r) \stackrel{\text{def}}{=} Y_s + Z_s(r) - r.$$

*If there exists an $r^* > 0$ such that $p(r^*) < 0$ then there exists a continuous curve $\{x(s)\}_{s \in [0,1]} \subset X$ such that $x(s) \in B_{r^*}(\hat{x}_s)$ and $H_s(x(s)) = 0$ for all $s \in [0,1]$. Furthermore, if $\hat{x}_\mathbf{s}, q_\mathbf{s} \in \mathcal{S}(X)$, then $x(s) \in \mathcal{S}(X)$.*

*Proof.* The first part of this Theorem has been presented and proven in many contexts, such as in [5] and [10].

The symmetry argument follows from the following considerations. Thanks to the injectivity of $A_s$, $T_s(x(s)) = x(s)$ implies $H_s(x(s)) = 0$. We have $H_s(x(s))^* = H_s(x(s)^*) = 0$ and thus $T_s(x(s)^*) = x(s)^*$. Since $\hat{x}_s \in \mathcal{S}(X)$, the ball $B_{r*}(\hat{x})$ is invariant under conjugate symmetry, hence $x(s)^* \in B_{r*}(\hat{x})$. Due to uniqueness of the fixed point of $T_s$ in $B_{r*}(\hat{x})$, we conclude that $x(s)^* = x(s) \in \mathcal{S}(X)$. □

When the continuation equation is chosen as in Remark 2.7, then the solution curve $\{x(s)\}_{s \in [0,1]} \subset X$ is smooth.

## 4 Linear operators

### 4.1 Operators on $X$

For simplicity, we introduce the notation $X = X_1 \times X_2 \times X_3$, where $X_1 = \mathbb{C}^3$, $X_2 = \ell^1_{\nu_2}$, $X_3 = \mathcal{L}^1_\nu$. Let $B \in B(X,X)$,

$$B : X \to X$$
$$x = (a,b,c) \mapsto Bx = (a',b',c'),$$

which we write as

$$\begin{pmatrix} a' \\ b' \\ c' \end{pmatrix} = \begin{pmatrix} B_{11} & B_{12} & B_{13} \\ B_{21} & B_{22} & B_{23} \\ B_{31} & B_{32} & B_{33} \end{pmatrix} \begin{pmatrix} a \\ b \\ c \end{pmatrix} \tag{31}$$

with

$$B_{ij} \in B(X_i, X_j).$$

We denote by $\|\cdot\|_{ij}$ the operator norm induced by the norms on $X_i$ and $X_j$. In particular, all the operator norms here introduced have a simple explicit form

or bound:

$$
\begin{aligned}
\|B_{11}\|_{11} &= \max_{n=1,2,3} \sum_{m=1,2,3} |B_{nm}|, \\
\|B_{12}\|_{12} &= \max_{n=1,2,3} \sup_{k\in\mathbb{Z}} \nu^{-|k|} |B_{nk}|, \\
\|B_{13}\|_{13} &= \max_{n=1,2,3} \sup_{(j,k)\in\mathbb{Z}^2} \nu^{-|k|-|j|} |B_{nk}|, \\
\|B_{21}\|_{21} &\leq \sum_{m=1,2,3} \sum_{k\in\mathbb{Z}} \nu^{|k|} |B_{km}|, \\
\|B_{22}\|_{22} &= \sup_{k_1\in\mathbb{Z}} \nu_2^{-|k_1|} \sum_{k\in\mathbb{Z}} \nu_2^{|k|} |B_{kk_1}| \\
\|B_{31}\|_{31} &\leq \sum_{m=1,2,3} \sum_{(j,k)\in\mathbb{Z}^2} \nu^{|k|+|j|} |B_{jkm}| \\
\|B_{32}\|_{32} &= \sup_{k_1\in\mathbb{Z}} \nu_2^{-|k_1|} \sum_{j,k\in\mathbb{Z}} \nu_1^{|j|}\nu_2^{|k|} |B_{jkk_1}| \\
\|B_{33}\|_{33} &= \sup_{j_1,k_1\in\mathbb{Z}} \nu_1^{-|j_1|}\nu_2^{-|k_1|} \sum_{j,k\in\mathbb{Z}} \nu_1^{|j|}\nu_2^{|k|} |B_{jj_1kk_1}|.
\end{aligned}
\tag{32}
$$

The norm $\|\cdot\|_{23}$ is not presented because it is never used in this paper.

On $B(X, X')$, we introduce the component norm

$$
\|B\|_c = \begin{pmatrix} \|B_{11}\|_{11} & \|B_{12}\|_{12} & \|B_{13}\|_{13} \\ \|B_{21}\|_{21} & \|B_{22}\|_{22} & \|B_{23}\|_{23} \\ \|B_{31}\|_{31} & \|B_{32}\|_{32} & \|B_{33}\|_{33} \end{pmatrix} \in \mathbb{R}^{3\times 3}.
$$

By the definition of operator norm on the space $X$, we have

$$
\|B\| \leq \max_{j=1,\ldots,3} \sum_{i=1,\ldots,3} \|B_{ij}\|_{ij}. \tag{33}
$$

## 4.2 Convolution operators on $\ell^1_{\nu_2}$ and $\mathcal{L}^1_\nu$

For $u = (u_k)_{k\in\mathbb{Z}} \in \ell^1_{\nu_2}$, we introduce the convolution operator as the operator $C^{\ell^1_{\nu_2}} : \ell^1_{\nu_2} \to B(\ell^1_{\nu_2}, \ell^1_{\nu_2})$ that satisfies

$$
C^{\ell^1_{\nu_2}}(u)w = u * w \qquad \text{for all } u, w \in \ell^1_{\nu_2}. \tag{34}
$$

By the definition of convolution, its matrix representation is

$$C^{\ell^1_{\nu_2}}(u)_{k_1 k_2} = u_{k_1-k_2}.$$

It follows from the Banach algebra property that

$$\|C^{\ell^1_{\nu_2}}(u)\|_{B(\ell^1_{\nu_2},\ell^1_{\nu_2})} = \|u\|_{\ell^1_{\nu_2}}. \tag{35}$$

Following the same structure in $\mathcal{L}^1_\nu$, we introduce the convolution operator $C^{\mathcal{L}^1_\nu} : \mathcal{L}^1_\nu \to B(\mathcal{L}^1_\nu, \mathcal{L}^1_\nu)$ such that

$$C^{\mathcal{L}^1_\nu}(u)w = u * w \qquad \text{for all } u, w \in \mathcal{L}^1_\nu. \tag{36}$$

We now have $C^{\mathcal{L}^1_\nu}(u)_{j_1 j_2 k_1 k_2} = u_{j_1-j_2, k_1-k_2}$, and

$$\|C^{\mathcal{L}^1_\nu}(u)\|_{B(\mathcal{L}^1_\nu, \mathcal{L}^1_\nu)} = \|u\|_{\mathcal{L}^1_\nu}. \tag{37}$$

## 4.3 Finite dimensional projection

In order to study a finite dimensional subset of $X$, we define for any $K = (M, N) \in \mathbb{N}^2$ the finite set of indices smaller than $K$ as

$$\mathbf{F}_K = \mathbf{F}_{(M,N)} = \{(j,k) \in \mathbb{Z}^2 \mid |j| \leq M, |k| \leq N\}. \tag{38}$$

The projection on the first $K = (M, N) \in \mathbb{N}^2$ modes is defined as

$$\begin{aligned}\Pi^{(K)} : \quad &X = \mathbb{C}^3 \times \ell^1_{\nu_2} \times \mathcal{L}^1_\nu \to X \\ &x = (a, b, c) \mapsto \Pi^{(K)}x = (a, \Pi^{(K)}_1 b, \Pi^{(K)}_2 c)\end{aligned} \tag{39}$$

with

$$(\Pi^{(K)}_1 b)_k = \begin{cases} b_k & \text{if } |k| \leq N, \\ 0 & \text{if } |k| > N, \end{cases} \tag{40}$$

and

$$(\Pi^{(K)}_2 c)_{jk} = \begin{cases} c_{jk} & \text{if } (j,k) \in \mathbf{F}_{(M,N)}, \\ 0 & \text{if } (j,k) \notin \mathbf{F}_{(M,N)}. \end{cases} \tag{41}$$

If $x = \Pi^{(K)}x \in X$, we say that $x$ has at most $K$ non-zero modes, and similarly for $b = \Pi_1^{(K)}b \in \ell^1_{\nu_2}$ and $c = \Pi_2^{(K)}c \in \mathcal{L}^1_\nu$. The complementary projection on the infinite tail is given by

$$\Pi^\infty_{(K)}x = x - \Pi^{(K)}x, \tag{42}$$

where the notation is meant to indicate the projection from the $(K+1)$-th mode onwards. There is a natural correspondence between $\Pi_1^{(K)}\ell^1_{\nu_2}$ and $\mathbb{C}^{2N+1}$, which we denote by the truncation operator

$$\tau_1^{(K)} : \ell^1_{\nu_2} \to \mathbb{C}^{2N+1}, \tag{43}$$
$$b \mapsto (b_k)_{k=-N,\ldots,N}, \tag{44}$$

as well as between $\Pi_2^{(K)}\mathcal{L}^1_\nu$ and $\mathbb{C}^{(2M+1)\times(2N+1)}$, which we indicate by

$$\tau_2^{(K)} : \mathcal{L}^1_\nu \to \mathbb{C}^{(2M+1)\times(2N+1)}, \tag{45}$$
$$c \mapsto (c_{jk})_{(j,k)\in \mathbf{F}_K}. \tag{46}$$

We introduce the notation for truncation in $X$ as

$$\tau^{(K)} : \quad X \to \mathbb{C}^3 \times \mathbb{C}^{2N+1} \times \mathbb{C}^{(2M+1)\times(2N+1)} \tag{47}$$
$$x = (a,b,c) \mapsto \tau^{(K)}x = (a, \tau_1^{(K)}b, \tau_2^{(K)}c). \tag{48}$$

Its corresponding inverse embedding is denoted by

$$E : \mathbb{C}^3 \times \mathbb{C}^{2N+1} \times \mathbb{C}^{(2M+1)\times(2N+1)} \to \Pi^{(K)}X.$$

**Definition 4.1.** $A \in \mathcal{B}(X,X)$ *is a* finite linear operator *of center dimensions $K_1$ and $K_2$ if*

$$Ax = \Pi^{(K_1)}(A\Pi^{(K_2)}x), \quad \text{for all } x \in X.$$

*If $K_1 = K_2 = K$, we say that $A$ has center dimension $K$.*

**Remark 4.2.** *A finite linear operator $A \in B(X,X)$ of center dimensions $K_1 = (M_1, N_1)$ and $K_2 = (M_2, N_2)$ can be represented by a complex matrix of dimension $(3+(2N_2+1)+(2M_2+1)(2N_2+1)) \times (3+(2N_1+1)+(2M_1+1)(2N_1+1))$.*

A class of operators we will be working with extensively is the class of eventually diagonal operators in $X$. It is an extention of the definition provided in [5] for eventually diagonal operators in $\ell^1_\nu$.

**Definition 4.3.** Let $A : X \to X$ be a linear operator. $A$ is eventually diagonal of center dimension $K$ if it can be decomposed into $A^F$ and $A^I$ such that
$$Ax = A^F x + A^I \Pi_{(K)}^\infty y$$
with $A^I$ a diagonal operator and $A^F$ a finite linear operator of center dimension $K$. The same definition applies for $A : \ell_\nu^1 \to \ell_\nu^1$ and $A : \mathcal{L}_\nu^1 \to \mathcal{L}_\nu^1$.

**Lemma 4.4.** *If $A$ is eventually diagonal,*
$$\|A_{ij}\|_{ij} \leqslant \max\{\|A_{ij}^F\|_{ij}, \|A_{ij}^I\|_{ij}\}, \qquad \text{for } i,j = 1,2,3.$$

*Proof.* Let $A = A^I + A^F$ be an eventually diagonal operator of center dimension $K = (M,N)$. For $i \neq j$, $A_{ij}^I = 0$, and $A_{11}^I = 0$, hence in these cases $\|A_{ij}\| = \|A_{ij}^F\|$. We are left to compute $\|A_{22}\|_{22}$ and $\|A_{33}\|_{33}$. By the definition of the norms in (32), we have

$$\begin{aligned}
\|A_{22}\|_{22} &= \sup_{k_2 \in \mathbb{Z}} \nu_2^{-|k_2|} \sum_{k_1 \in \mathbb{Z}} \nu_2^{|k_1|} |(A_{22})_{k_1 k_2}| \\
&= \max\{ \sup_{|k_2| \leqslant N} \nu_2^{-|k_2|} \sum_{k_1 \in \mathbb{Z}} \nu_2^{|k_1|} |(A_{22})_{k_1 k_2}|, \\
&\qquad \sup_{|k_2| > N} \nu_2^{-|k_2|} \sum_{k_1 \in \mathbb{Z}} \nu_2^{|k_1|} |(A_{22})_{k_1 k_2}|\} \\
&= \max\{ \sup_{|k_2| \leqslant N} \nu_2^{-|k_2|} \sum_{k_1 \in \mathbb{Z}} \nu_2^{|k_1|} |(A_{22}^F)_{k_1 k_2}|, \sup_{|k_2| > N} |(A_{22}^I)_{k_1 k_1}|\} \\
&= \max\{\|A_{22}^F\|_{22}, \|A_{22}^I\|_{22}\}.
\end{aligned}$$

The same approach can be applied to $A_{33}$, thus concluding the proof. $\square$

**Remark 4.5.** *The product of two eventually diagonal operators of center dimensions $K_1$ and $K_2$, respectively, is an eventually diagonal operator of center dimension*
$$K = \max\{K_1, K_2\}.$$

The operator $A$ introduced in Section 3 will be constructed to be an eventually diagonal operator. We therefore want to study this type of operators a little more in depth.

**Remark 4.6.** Let $x \in X$ such that $x = \Pi^{(K')}x$, and let $A$ be an eventually diagonal operator of center dimension $K$, then
$$Ax = \Pi^{(\tilde{K})}(Ax),$$
with
$$\tilde{K} = \max(K, K').$$

Another case in which it is easy to compute the number of non-zero modes is when we are considering convolutions.

**Lemma 4.7.** Let $x = \Pi_2^{(K_1)}x$ and $v = \Pi_2^{(K_2)}v$ then
$$x * v = \Pi_2^{(K_1+K_2)}(x * v).$$
The same holds in $\ell^1_{\nu_2}$.

We also introduce the definition of a $K$-diagonal operator, where $K$ can be thought of as a band width.

**Definition 4.8.** $B : X \to X$ is called a $K$-diagonal operator on $X$ if
$$By = \Pi^{(K+K')}(By)$$
for all $y \in X$ such that $y = \Pi^{(K')}y$. And similarly for $B : \mathcal{L}^1_\nu \to \mathcal{L}^1_\nu$ and $B : \ell^1_{\nu_2} \to \ell^1_{\nu_2}$.

**Remark 4.9.** Diagonal operators are 1-diagonal operators. Furthermore, by Lemma 4.7, the convolution operators associated to $y \in \ell^1_{\nu_2}$ and $z \in \mathcal{L}^1_\nu$ with $K$ non-zero modes are $K$-diagonal.

**Lemma 4.10.** $B$ is a $K$-diagonal operator in $\mathcal{L}^1_\nu$ if and only if its matrix representation $B$ is such that $B_{jj_1kk_1} = 0$ if $(j-j_1, k-k_1) \notin \mathbf{F}_K$. A similar statement holds in $\ell^1_{\nu_2}$.

**Lemma 4.11.** Let $K = (M, N) \in \mathbb{N}^2$, and let $B : X \to X$ be a $K$-diagonal linear operator. Then $B_{22} : \ell^1_{\nu_2} \to \ell^1_{\nu_2}$ is a $K$-diagonal linear operator and so is $B_{33} : \mathcal{L}^1_\nu \to \mathcal{L}^1_\nu$. For $B_{32} : \ell^1_{\nu_2} \to \mathcal{L}^1_\nu$, its matrix representation $B_{32}$ is such that $(B_{32})_{jkk_1} = 0$ if $|k - k_1| > N$.

**Theorem 4.12.** Let $A : X \to X$ be a finite dimensional operator of center dimension $K_1 = (M_1, N_1)$ and $K_2 = (M_2, N_2)$. Let $B : X \to X$ be a $\tilde{K}$-diagonal operator. Then the operator product $AB$ is a finite dimensional operator of center dimension $K_1$ and $K_2 + \tilde{K}$.

*Proof.* The proof will be carried out for $A' = A_{22}$ and $B' = B_{22}$, but its extension to the whole of $X$ follows the same reasoning.

By definition of the operator product, it holds that

$$(A'B')_{k_1 k_2} = \sum_{k_3 \in \mathbb{Z}} A'_{k_1 k_3} B'_{k_3 k_2}.$$

The assumptions imply that $A'_{k_1 k_3} = 0$ if $k_1 > N_1$ or $k_3 > N_1$. Furthermore, $B'_{k_3 k_2} = 0$ if $|k_3 - k_2| > \tilde{N}$. Hence, we can rewrite

$$(A'B')_{k_1 k_2} = \sum_{k_3 \in I_{k_2}} A'_{k_1 k_3} B'_{k_3 k_2},$$

with

$$I_{k_2} = \{k_2 - \tilde{N}, \ldots, k_2 + \tilde{N}\} \cap \{-N_2, \ldots, N_2\}.$$

The set $I_{k_2}$ is non-empty only if $|k_2| \leq N_2 + \tilde{N}$. This concludes the proof. $\square$

**Corollary 4.13.** *Under the hypothesis of Theorem 4.12, we have*

$$AB = \Pi^{(K_1)} A \Pi^{(K_2)} B = \Pi^{(K_1)} A \Pi^{(K_2)} B \Pi^{(K_2 + \tilde{K})}.$$

*Hence, (a bound on) the operator norm $\|AB\|$ is explicitly computable.*

# 5 Construction of the operator $A_s$

In this Section, we construct the operator $A_s$ needed for the application of the radii polynomial approach to the Kuramoto-Sivashinsky PDE given by (24). We consider the segment

$$\hat{x}_s = E\hat{x}_0 + s(E\hat{x}_1 - E\hat{x}_0),$$

where $\hat{x}_0, \hat{x}_1 \in \tau^{(K)} \mathcal{S}(X)$ are the two numerical approximate solutions of (24) already introduced at the end of Section 2.

The radii polynomial approach to the validation of a branch relies strongly on finding a reasonably good approximation $A_s$ of the inverse of $DH_s(\hat{x}_s)$. To construct $A_s$, we will introduce the linear operator $A_s^\dagger$, an approximation of

$DH_s(\hat{x}_s)$. First, let us compute explicitly

$$DH_s(x) =$$
$$\tilde{\Pi} \begin{pmatrix} D_{(\lambda_1,\lambda_2,a)} \begin{pmatrix} \boldsymbol{E}_s(x) \\ G(x) \end{pmatrix} & D_y \begin{pmatrix} \boldsymbol{E}_s(x) \\ G(x) \end{pmatrix} & D_z \begin{pmatrix} \boldsymbol{E}_s(x) \\ G(x) \end{pmatrix} \\ D_{(\lambda_1,\lambda_2,a)} F_1(x) & \lambda_2 K^4 - K^2 - 2iKC^{\ell_{\nu_2}^1}(y) & 0 \\ D_{(\lambda_1,\lambda_2,a)} F_2(x) & -2\lambda_1 iKC^{\mathcal{L}_\nu^1}(z)E_\ell & \begin{matrix} iJ + \lambda_1\lambda_2 K^4 - \lambda_1 K^2 \\ -2\lambda_1 iKC^{\mathcal{L}_\nu^1}(E_\ell y + az) \end{matrix} \end{pmatrix}. \tag{49}$$

with

$$D_{(\lambda_1,\lambda_2,a)} F_1(x) = \begin{pmatrix} 0 & K^4 y & 0 \end{pmatrix},$$
$$D_{(\lambda_1,\lambda_2,a)} F_2(x) = \begin{pmatrix} (\lambda_2 K^4 - K^2)z & \lambda_1 K^4 z & -\lambda_1 iK(z * z) \\ -iK(2z * E_\ell y + a(z*z)) & & \end{pmatrix},$$

We define
$$\hat{H}_s(x) \stackrel{\text{def}}{=} \tau^{(K)} H_s(\tau^{(K)} x),$$

which can be interpreted as a map from $\mathbb{C}^{3+(2N+1)+(2N+1)(2M+1)}$ to itself, and we have
$$D\hat{H}_s(x) = \tau^{(K)} DH_s(\tau^{(K)} x)\tau^{(K)}, \tag{50}$$

which can be represented by a complex square matrix of dimensions $(3 + (2N + 1) + (2N + 1)(2M + 1))$.

We construct the dependence of $A_s^\dagger$ on $s$ to be linear:
$$A_s^\dagger = A_0^\dagger + s(A_1^\dagger - A_0^\dagger),$$

where we define $A_{\mathsf{s}}^\dagger$ as
$$A_{\mathsf{s}}^\dagger = ED\hat{H}_{\mathsf{s}}(\hat{x}_{\mathsf{s}}) + A_{\mathsf{s}}^{\dagger I}, \qquad \text{for } \mathsf{s} = 0, 1,$$

with $A_{\mathsf{s}}^{\dagger I}$ the tail operator
$$A_{\mathsf{s}}^{\dagger I} = \begin{pmatrix} 0 & 0 & 0 \\ 0 & \hat{\lambda}_{2\mathsf{s}} K^4 - K^2 & 0 \\ 0 & 0 & iJ + \hat{\lambda}_{1\mathsf{s}}\hat{\lambda}_{2\mathsf{s}} K^4 - \hat{\lambda}_{1\mathsf{s}} K^2 \end{pmatrix} \Pi_{(K)}^\infty.$$

For the approximate inverse $A_s$ of $A_s^\dagger$, we set
$$A_s = A_0 + s(A_1 - A_0),$$

and
$$A_{\mathsf{s}} = A_{\mathsf{s}}^F + A_{\mathsf{s}}^I, \qquad \text{for } \mathsf{s} = 0, 1,$$
with
$$A_{\mathsf{s}}^F = EQ_s\tau^{(K)}, \qquad \text{for } \mathsf{s} = 0, 1,$$
where $Q_s$ is an approximate numerical inverse of the matrix representation of $D\hat{H}_{\mathsf{s}}(\tau^{(K)}x_{\mathsf{s}})$, a complex square matrix of dimension $(3 + (2N + 1) + (2N + 1)(2M + 1))$. The diagonal operators are defined as

$$A_{\mathsf{s}}^I = \begin{pmatrix} 0 & 0 & 0 \\ 0 & A_{\mathsf{s}1}^I & 0 \\ 0 & 0 & A_{\mathsf{s}2}^I \end{pmatrix} \Pi_{(K)}^{\infty}, \tag{51}$$

with
$$A_{\mathsf{s}1}^I = (\hat{\lambda}_{2\mathsf{s}}\mathrm{K}^4 - \mathrm{K}^2)^{-1}$$
and
$$A_{\mathsf{s}2}^I = (\mathrm{i}\mathrm{J} + \hat{\lambda}_{1\mathsf{s}}\hat{\lambda}_{2\mathsf{s}}\mathrm{K}^4 - \hat{\lambda}_{1\mathsf{s}}\mathrm{K}^2)^{-1},$$
where the inverse here is interpreted elementwise, that is, for $(j, k) \notin \mathbf{F}_K$,
$$(A_{\mathsf{s}1}^I y)_k = \frac{y_k}{\hat{\lambda}_{2\mathsf{s}}k^4 - k^2},$$
and
$$(A_{\mathsf{s}2}^I z)_{jk} = \frac{z_{jk}}{\mathrm{i}j + \hat{\lambda}_{1\mathsf{s}}\hat{\lambda}_{2\mathsf{s}}k^4 - \hat{\lambda}_{1\mathsf{s}}k^2}.$$

In Equation (28) we stated two conditions which $A_s$ needs to satisfy. First, we require that $A_s H_s(x) \in X$ for all $x \in X$. Although $H_s(x)$ is not in $X$ due to the unboundedness of the operators J and K, this is counteracted by the action of the tail of $A_s$, hence $A_s H_s(x) \in X$.

The second condition is that $A_s$ is injective. This condition is proven *a posteriori* by the radii polynomial approach. Indeed, if the validation succeeds, it holds that $1 > Z(r) \geq \sup_{s \in [0,1]} \|I - A_s A_s^\dagger\|$, thus $A_s$ is surjective for all $s \in [0, 1]$. We also know that $A_s$ is eventually diagonal for all $s$, that is $A_s = A_s^F + A_s^I$. We now interpret this as a block-diagonal decomposition of $A_s$. Then, by surjectivity of $A_s$, we deduce surjectivity of $A_s^F$, and since $A_s^F$ is finite dimensional, surjectivity implies injectivity. Injectivity of $A_s^I$ follows from the nontriviality of the elements on its diagonal. Hence $A_s$ is injective.

**Remark 5.1.** *The elements on the diagonal are nontrivial if $\hat{\lambda}_{2\mathsf{s}}k^4 - k^2 \neq 0$ for all $|k| > N$. Thus we impose*
$$N^2 > \max\{\hat{\lambda}_{20}^{-1}, \hat{\lambda}_{21}^{-1}\}.$$

# 6 Bounds for the norm of eventually diagonal operators

This section is devoted to bound the norm of some eventually diagonal operators, where we add some constraint on the diagonal part of the operator. Let $P$ be an eventually diagonal operator in $B(\mathcal{L}_\nu^1, \mathcal{L}_\nu^1)$ of center dimension $K = (M, N) \in \mathbb{N}^2$, then it has the form

$$P : \mathcal{L}_\nu^1 \to \mathcal{L}_\nu^1 \tag{52}$$
$$u \mapsto P : (Pu)_{jk} = \begin{cases} E(B\tau_2^K u)_{jk} & \text{for } (j,k) \in \mathbf{F}_K, \\ \mathbf{p}(j,k)u_{jk} & \text{otherwise,} \end{cases}$$

for a function $\mathbf{p}$ and a matrix $B \in \mathbb{C}^{(2M+1)(2N+1) \times (2M+1)(2N+1)}$, see Remark 4.2. In this section we will study the norm of $P$ when $\mathbf{p}$ is a rational polynomial of the form

$$\mathbf{p}(j,k) = \frac{q(j,k)}{p(j,k)} = \frac{q_1(j) + q_2(k)}{\mathrm{i}p_1(j) + p_2(k)}, \tag{53}$$

with $p_1(j)$ and $p_2(k)$ real-valued polynomials and $q_1(j)$ and $q_2(k)$ complex polynomials. If there is any constant term in $q$, it will be assigned to $q_2$.

It holds, by Lemma 4.4 and by the definitions of the operator norms in (32), that

$$\|P\| = \max\{\|B\|, \sup_{(j,k) \notin \mathbf{F}_K} |\mathbf{p}(j,k)|\}.$$

**Remark 6.1.** *The bound on $\sup_{(j,k) \notin \mathbf{F}_K} |\mathbf{p}(j,k)|$ can be finite only if the order of $q_i$ is smaller than or equal to the order of $p_i$, for $i = 1, 2$.*

To simplify the presentation, we will further restrict ourselves to $p_1(j) = j$, although the results can be generalised to polynomial $p_1(j)$.

By the triangle inequality

$$\sup_{(j,k) \notin \mathbf{F}_K} |\mathbf{p}(j,k)| \leq \sup_{(j,k) \notin \mathbf{F}_K} \left| \frac{q_1(j)}{\mathrm{i}j + p_2(k)} \right| + \sup_{(j,k) \notin \mathbf{F}_K} \left| \frac{q_2(k)}{\mathrm{i}j + p_2(k)} \right| \tag{54}$$

We will consider the two fractions separately. With the restrictions imposed on $q_1$, we have $q_1(j) = c_1 j$. Therefore the first fraction is bounded trivially by $|c_1|$.

**Remark 6.2.** *If there exists a $(j^*, k^*) \notin F_K$ such that $p(j^*, k^*) = 0$, then the norm of $P$ is infinite. Hence in the following we want to ensure that there is no $k^* \in \mathbb{N}$, such that $|k^*| > N$ and $p(0, k^*) = 0$. This is incorporated in the definition of $Q_N$ in (59) below, which will simply be infinity if this is violated.*

Concerning the second term in the righthand side of (54), we remark that, for fixed $k$, $\left|\frac{q_2(k)}{ij+p(k)}\right|$ decreases if $|j|$ increases, hence

$$\sup_{(j,k)\notin \mathbf{F}_K} \left|\frac{q_2(k)}{ij + p_2(k)}\right| \leqslant \max\left\{\max_{|j|=M+1, |k|\leqslant N+1}\left|\frac{q_2(k)}{ij+p_2(k)}\right|, \sup_{|k|>N}\left|\frac{q_2(k)}{p_2(k)}\right|\right\}. \tag{55}$$

The first maximum of (55) is taken over a finite set, hence it is directly computable. The supremum, on the other hand, is computed over an infinite set and requires a more refined approach, i.e., we want to bound

$$\sup_{|k|>N}\left|\frac{q_2(k)}{p_2(k)}\right| \tag{56}$$

using a finite number of computations. As a first step, we construct $k_1' \in \mathbb{N}$ satisfying

$$p_2(x) \neq 0 \quad \text{for all } |x| > k_1', \tag{57}$$

and $k_2' \in \mathbb{N}$ satisfying

$$q_2'(x)p_2(x) - p_2'(x)q_2(x) = 0 \quad \text{for all } |x| > k_2'. \tag{58}$$

To find $k_1$ and $k_2$, we use the minimum of Lagrange's and Cauchy's bounds, which provide explicit finite bounds on roots of polynomials. Then we set $k' = \max\{k_1', k_2'\}$.

It then holds that (56) is bounded by

$$Q_N \stackrel{\text{def}}{=} \max\left\{\lim_{k\to\pm\infty}\left|\frac{q_2(k)}{p_2(k)}\right|, \max_{|k|>N, |k|\leqslant k'}\left|\frac{q_2(k)}{p_2(k)}\right|, \left|\frac{q_2(N+1)}{p_2(N+1)}\right|, \left|\frac{q_2(-N-1)}{p_2(-N-1)}\right|\right\}. \tag{59}$$

The limit $\lim_{k\to\pm\infty}\left|\frac{q_2(k)}{p_2(k)}\right|$ is trivial to determine and is non-zero only if $q_2$ and $p_2$ have the same polynomial order. All other bounds only require a finite computation.

With this algorithm, we can then bound explicitly the norm of an eventually diagonal operator as defined in (52), (53).

# 7 Single solution bounds for radii polynomial

For ease of exposition, we will first treat the case of a single solution. This may also be interpreted as finding the bounds required for the **s**-independent version of Theorem 3.1, i.e., the case where $\hat{x}_0 = \hat{x}_1$ and $\hat{q}_0 = \hat{q}_1$. We note that although the arguments below are self-contained, some additional details can be found in [9].

## 7.1 $Y$ bound

The $Y$ bound is an upper bound to
$$\|T(\hat{x}) - \hat{x}\|_X = \|AH(\hat{x})\|_X.$$

By construction, $\hat{x}$ has $K$ non-zero modes. Furthermore, $F_1$ and $F_2$ are the sum of linear terms and second order convolutions. Thus $F_i(\hat{x}) = \Pi_i^{(2K)} F_i(\hat{x})$ for $i = 1, 2$, and as a consequence $H(\hat{x}) = \Pi^{(2K)} H(\hat{x})$. We know that $A$ is eventually diagonal of center dimension $K$ and Remark 4.6 can be applied. Therefore $AH(\hat{x})$ has $2K$ non-zero modes. Thus, its norm can be computed directly with a finite number of computations in interval arithmetic.

## 7.2 $Z$ bound

The $Z$ bound presented in (30) can be split into three parts by considering
$$DT(\hat{x} + rb)rc = \left[(I - AA^\dagger) + A(A^\dagger - DH(\hat{x})) + A(DH(\hat{x}) - DH(\hat{x} + rb))\right] rc,$$
and computing their bounds separately. We write
$$Z(r) = Z_0 r + Z_1 r + Z_2(r) r \tag{60}$$
with
$$Z_0 \geqslant \|I - AA^\dagger\|_{B(X,X)}, \tag{61}$$
$$Z_1 \geqslant \|A(A^\dagger - DH(\hat{x}))\|_{B(X,X)}, \tag{62}$$
and
$$Z_2(r) \geqslant \|A(DH(\hat{x}) - DH(\hat{x} + rb))\|_{B(X,X)}, \tag{63}$$
where $I \in B(X, X)$ is the identity operator.

## 7.3 $Z_0$ bound

In the case at hand, the $Z_0$ bound is fairly simple, because $A$ and $A^\dagger$ are eventually diagonal operators of center dimension $K$. Thus, applying the triangle inequality we have

$$\|I - AA^\dagger\|_{B(X,X)} = \|I\Pi_{(K)}^\infty - A^I A^{\dagger I}\Pi_{(K)}^\infty\|_{B(X,X)} + \|I\Pi^{(K)} - A^F A^{\dagger F}\|_{B(X,X)}.$$

We remark that, by definition, the tails of the operator $A$ and $A^\dagger$ are exact inverses of one another. Therefore, the tail of $A^I A^{\dagger I}$ is the identity, and the first norm is zero.

We are left to compute

$$\|I\Pi^{(K)} - A^F A^{\dagger F}\|_{B(X,X)}.$$

The matrix representations of $A^F$ and $A^{\dagger F}$ are known and finite. Hence we can compute this norm applying Equation (33) and the definition of the operator norms in (32).

## 7.4 $Z_1$ bound

In the computation of the $Z_1$ bound, the subtraction $DH(\hat{x}) - A^\dagger$ gives an operator where most terms cancel and which acts as a convolution. We can represent it as

$$\begin{aligned}
B &= DH(\hat{x}) - A^\dagger \\
&= \begin{pmatrix} 0 & 0 & 0 \\ 0 & \mathrm{K}b_{22} & 0 \\ \mathrm{K}b_{31} & \mathrm{K}b_{32} & \mathrm{K}b_{33} \end{pmatrix} - \Pi^{(K)} \begin{pmatrix} 0 & 0 & 0 \\ 0 & \mathrm{K}b_{22} & 0 \\ \mathrm{K}b_{31} & \mathrm{K}b_{32} & \mathrm{K}b_{33} \end{pmatrix} \Pi^{(K)},
\end{aligned} \qquad (64)$$

with

$$\begin{aligned}
b_{22} &= 2\mathrm{i}C^{\ell^1_{\nu_2}}(\hat{y}), \\
b_{31} &= \left[\mathrm{i}(2\hat{z} * E_\ell \hat{y} + \hat{a}(\hat{z}*\hat{z})), 0, \hat{\lambda}_1 \mathrm{i}(\hat{z}*\hat{z})\right], \\
b_{32} &= 2\hat{\lambda}_1 \mathrm{i} C^{\mathcal{L}^1_\nu}(\hat{z}) E_\ell, \\
b_{33} &= 2\hat{\lambda}_1 \mathrm{i} C^{\mathcal{L}^1_\nu}(E_\ell \hat{y} + \hat{a}\hat{z}),
\end{aligned}$$

where, with some abuse of notation,

$$\mathrm{K}b_{31} = \mathrm{K}\left[\mathrm{i}(2\hat{z} * E_\ell \hat{y} + \hat{a}(\hat{z} * \hat{z})), 0, \hat{\lambda}_1 \mathrm{i}(\hat{z} * \hat{z})\right]$$
$$= \left[\mathrm{Ki}(2\hat{z} * E_\ell \hat{y} + \hat{a}(\hat{z} * \hat{z})), 0, \mathrm{K}\hat{\lambda}_1 \mathrm{i}(\hat{z} * \hat{z})\right].$$

**Remark 7.1.** *In the second block-row of B the row corresponding to $k = 0$, which was âĂIJredefinedâĂİ through $\tilde{\Pi}$ to correspond to $G_3^{\hat{x}}$ rather than $(F_1)_0$, vanishes. Similar considerations hold for the the row corresponding to $j = k = 0$ in the third block-row of B. Hence the compact notation in (7.5) may be used without including the projection $\tilde{\Pi}$. In the same way, in the rest of the paper, when the rows corresponding to $k = 0$ and $j = k = 0$ vanish, we drop the use of the projection $\tilde{\Pi}$, such as in Equation (65) as well as for higher derivatives of H.*

We decompose

$$A(A^\dagger - DH(\hat{x})) = A^F(A^\dagger - DH(\hat{x})) + A^I(A^\dagger - DH(\hat{x})),$$

hence

$$\|A(A^\dagger - DH(\hat{x}))\| \leqslant \|A^F(A^\dagger - DH(\hat{x}))\| + \|A^I(A^\dagger - DH(\hat{x}))\|.$$

The first term has a finite number of non-zero elements, by Theorem 4.12. The norm can therefore be computed explicitly. The second term can be handled by computing explicitly the product $A^I B$ with the help of the two notations introduced in (51) and (64). Then we have

$$A^I B = \begin{pmatrix} 0 & 0 & 0 \\ 0 & A_1^I \mathrm{K} b_{22} & 0 \\ A_2^I \mathrm{K} b_{31} & A_2^I \mathrm{K} b_{32} & A_2^I \mathrm{K} b_{33} \end{pmatrix} - \Pi^{(K)} \begin{pmatrix} 0 & 0 & 0 \\ 0 & A_1^I \mathrm{K} b_{22} & 0 \\ A_2^I \mathrm{K} b_{31} & A_2^I \mathrm{K} b_{32} & A_2^I \mathrm{K} b_{33} \end{pmatrix} \Pi^{(K)},$$

where the second term vanishes, and it holds that

$$\|A^I B\|_c \leqslant \begin{pmatrix} 0 & 0 & 0 \\ 0 & \|A_1^I \mathrm{K}\|_{22} \|b_{22}\|_{22} & 0 \\ \|A_2^I \mathrm{K}\|_{33} \|b_{31}\|_{31} & \|A_2^I \mathrm{K}\|_{33} \|b_{32}\|_{32} & \|A_2^I \mathrm{K}\|_{33} \|b_{33}\|_{33} \end{pmatrix},$$

with the inequality interpreted elementwise. By using Equations (35) and (37),

we obtain

$$\|b_{22}\|_{22} \leqslant 2\|\hat{y}\|_{\ell^1_{\nu_2}},$$
$$\|b_{31}\|_{31} \leqslant \max\{\|2\hat{z} * E_\ell \hat{y} + \hat{a}(\hat{z} * \hat{z})\|_{\mathcal{L}^1_\nu}, 0, |\hat{\lambda}_1|\|\hat{z} * \hat{z}\|_{\mathcal{L}^1_\nu}\},$$
$$\|b_{32}\|_{32} \leqslant 2|\hat{\lambda}_1|\|\hat{z}\|_{\mathcal{L}^1_\nu},$$
$$\|b_{33}\|_{33} \leqslant 2|\hat{\lambda}_1|\|E_\ell \hat{y} + \hat{a}\hat{z}\|_{\mathcal{L}^1_\nu},$$

and $\|A^I_i K\|$ for $i = 0, 1$ fits the setup discussed in Section 6 and can be computed via (59).

## 7.5  $Z_2$ bound

Following the approach already presented in [10], for the computation of $Z_2(r)$ we can apply the mean value theorem. We introduce $R$, an *a priori* bound on the validation radius $r^*$, and we obtain

$$\|A(DH(\hat{x}) - DH(\hat{x} + rb))\|_{B(X,X)} \leqslant \sup_{b,c \in B_1(0)} \|AD^2 H(\hat{x} + Rb)c\|_{B(X,X)}.$$

The major difference with respect to [10] is that in the case at hand the norm of $D^2 H(\hat{x} + Rb)$ is not finite because it includes spatial derivative operators. Therefore, in this case the multiplication with $A$ plays a key role in keeping the $Z_2(r)$ bound finite. We start by computing $D^2 H$ explicitly, using the notation introduced in (31). We compute it at a point $x = (\lambda_1, \lambda_2, a, y, z)$ and apply it to $w = (\beta_1, \beta_2, b, r, s) \in X$:

$$D^2 H(x)w = \begin{pmatrix} 0 & 0 & 0 \\ d_{21}(x,w) & d_{22}(x,w) & 0 \\ d_{31}(x,w) & d_{32}(x,w) & d_{33}(x,w) \end{pmatrix}, \tag{65}$$

with

$$d_{21}(x,w) = [0, K^4 r, 0],$$
$$d_{22}(x,w) = K^4 \beta_2 - 2iKC^{\ell^1_{\nu_2}}(r),$$
$$d_{31}(x,w) = d_{311}(x,w), K^4(z\beta_1 + \lambda_1 s), -iK(\beta_1 z * z + 2\lambda_1 z * s),$$
$$d_{311}(x,w) = K^4(z\beta_2 + \lambda_2 s) - K^2 s - iK(2E_\ell y * s + 2z * E_\ell r + bz * z + 2az * s),$$
$$d_{32}(x,w) = -2iKC^{\mathcal{L}^1_\nu}(\lambda_1 s + \beta_1 z) E_\ell,$$
$$d_{33}(x,w) = K^4(\lambda_2 \beta_1 + \lambda_1 \beta_2) - K^2 \beta_1$$
$$\qquad - 2iKC^{\mathcal{L}^1_\nu}(\beta_1 E_\ell y + \lambda_1 E_\ell r + a\beta_1 z + \lambda_1(bz + as)).$$

We notice that $D^2H$ is a $2K$-diagonal operator since it is the product of diagonal operators and convolution operators. All $d_{nm}(x,w)$ are of the form

$$d_{nm}(x,w) = \mathrm{K} d^{[1]}_{nm}(x,w) + \mathrm{K}^2 d^{[2]}_{nm}(x,w) + \mathrm{K}^4 d^{[4]}_{nm}(x,w), \tag{66}$$

for $n, m = 1, \ldots, 3$.

We use the notation introduced in (31) to write

$$A = \begin{pmatrix} A_{11} & A_{12} & A_{13} \\ A_{21} & A_{22} & A_{23} \\ A_{31} & A_{32} & A_{33} \end{pmatrix}.$$

With this notation, we can write

$$(AD^2H(\hat{x}+Rb)c)_{im} = \sum_{n=1,\ldots,3} A_{in} d_{nm}(\hat{x}+Rb,c),$$

and by the triangle inequality we can bound

$$\|(AD^2H(\hat{x}+Rb)c)_{im}\|_{im} \leqslant \sum_{n=1,\ldots,3} \|A_{in} d_{nm}(\hat{x}+Rb,c)\|_{im}.$$

With the notation introduced in (66) and applying again the triangle inequality, we have

$$\|A_{in} d_{nm}(\hat{x}+Rb,c)\|_{im} \leqslant \sum_{p=1,2,4} \|A_{in} \mathrm{K}^p\|_{in} \|d^{[p]}_{nm}(\hat{x}+Rb,c)\|_{nm}.$$

Since $K$ is a diagonal operator, it holds that, by Lemma 4.4,

$$\|A_{ij} \mathrm{K}^p\|_{ij} = \max\{\|A_{ij}^F \Pi^{(K)} \mathrm{K}^p\|_{ij}, \|A_{ij}^I \Pi^{\infty}_{(K)} \mathrm{K}^p\|_{ij}\},$$

where the first norm can be computed directly and the second can be estimated by applying the bounds from Section 6. The norms $\|d^{[p]}_{nm}(\hat{x}+Rb,c)\|_{nm}$ can be easily estimated for all $n$ and $m$ because they are sums and products of scalars and convolution operators, hence (35) and (37) and the bound $\|x+Rb\|_c \leqslant \|x\|_c + R$ suffice to bound these terms.

As an example we consider $d^{[p]}_{33}(x,w)$, which can be written as

$$d^{[1]}_{33}(x,w) = -2\mathrm{i}C^{\mathcal{L}^1_\nu}(\beta_1 E_\ell y + \lambda_1 E_\ell r + a\beta_1 z + \lambda_1(bz+as)),$$
$$d^{[2]}_{33}(x,w) = -\beta_1,$$
$$d^{[4]}_{33}(x,w) = \lambda_2 \beta_1 + \lambda_1 \beta_2.$$

Using $\|x + Rb\|_c \leqslant \|x\|_c + R$ and $\|b\|_c \leqslant 1$, we obtain, for any $\|w\|_X \leqslant 1$,

$$\|d_{33}^{[1]}(x + Rb, w)\|_{33} \leqslant 2\left(\|y\|_\ell + |\lambda_1| + 2R + (|a| + R)(\|z\|_{\mathcal{L}_\nu^1} + R)\right)$$
$$+ 2(|\lambda_1| + R)(\|z\|_{\mathcal{L}_\nu^1} + |a| + 2R),$$
$$\|d_{33}^{[2]}(x + Rb, w)\|_{33} \leqslant 1,$$
$$\|d_{33}^{[4]}(x + Rb, w)\|_{33} \leqslant |\lambda_2| + |\lambda_1| + 2R.$$

The estimates of the other terms are treated in an analogous manner.

## 8 Bounds for segment validation

To compute the bounds introduced in Section 3 that include a maximum over $s \in [0, 1]$, we use the following lemma.

**Lemma 8.1.** *Let $f : [0, 1] \mapsto \mathbb{C}$ be $C^2$, then*

$$\max_{s \in [0,1]} |f(s)| \leqslant \max\{|f(0)|, |f(1)|\} + \frac{1}{8} \max_{x \in [0,1]} |f''(x)|. \tag{67}$$

For the proof, we refer to [10, Lemma 8.1]. We want to apply this lemma to a variety of cases. The most general application is the following.

**Lemma 8.2.** *Let $i, j \in \{1, 2, 3\}$. For $s \in [0, 1]$ let $A(s) : X_i \to X_j$ be a family of bounded linear operators dependent on a parameter $s$, with matrix representation $A(s) = (A_{k_1 k_2}(s))_{k_1 \in I_i, k_2 \in I_j}$, with $I_1$ and $I_2$ the appropriate set of indices for $X_i$ and $X_j$ respectively, that is $I_1 = \{1, 2, 3\}, I_2 = \mathbb{Z}, I_3 = \mathbb{Z}^2$. Then*

$$\max_{s \in [0,1]} \|A(s)\|_{ij} \leqslant \left\|(\max\{|A_{k_1 k_2}(0)|, |A_{k_1 k_2}(1)|\})_{k_1 \in I_i, k_2 \in I_j}\right\|_{ij} +$$
$$\frac{1}{8} \left\|\left(\sup_{x \in [0,1]} |A''_{k_1 k_2}(x)|\right)_{k_1 \in I_i, k_2 \in I_j}\right\|_{ij},$$

*provided the elementwise second derivatives $A''_{k_1 k_2}$ exists for all $k_1 \in I_1$ and $k_2 \in I_2$ and the second term is finite.*

For the proof, we refer to [10]. That proof is carried out on a different space but the proof on $X$ follows in the same way.

Considering the space $X$ and the linear operator $B(s) : X \to X$ for $s \in [0, 1]$, using the notation (31), we have

$$\max_{s\in[0,1]} \|B(s)\|_{B(X,X)} \leqslant \max_{j=1,\ldots,3} \sum_{i=1,\ldots,3} \max_{s\in[0,1]} \|B_{ij}(s)\|, \qquad (68)$$

and we can apply Lemma 8.2 to each $B_{ij}$. Having already introduced the operators $A_0$ and $A_1$, we introduce a notation for their difference:

$$A_\Delta = A_1 - A_0.$$

In general, the subscript $\Delta$ will indicate the difference between the element at hand at $s = 1$ and $s = 0$.

**Remark 8.3.** *By construction, $A_\Delta$ is an eventually diagonal operator of center dimension $K$, hence $A_\Delta = A_\Delta^F + A_\Delta^I$. The tail of $A_\Delta$ is non-zero in general, but applying the same notation as in (51), it is computable element wise by*

$$A_\Delta^I = \begin{pmatrix} 0 & 0 & 0 \\ 0 & A_{\Delta 1}^I & 0 \\ 0 & 0 & A_{\Delta 2}^I \end{pmatrix} \Pi_{(K)}^\infty. \qquad (69)$$

*Here*

$$A_{\Delta 1}^I = P_1(1)^{-1} - P_1(0)^{-1} = \frac{P_1(0) - P_1(1)}{P_1(0)P_1(1)},$$

*and*

$$A_{\Delta 2}^I = P_2(1)^{-1} - P_2(0)^{-1} = \frac{P_2(0) - P_2(1)}{P_2(0)P_2(1)},$$

*with*

$$P_1(\mathsf{s}) \stackrel{\text{def}}{=} \hat{\lambda}_{\mathsf{s}2}\mathrm{K}^4 - \mathrm{K}^2, \quad \textit{for } \mathsf{s} = 0, 1,$$

*and*

$$P_2(\mathsf{s}) \stackrel{\text{def}}{=} \mathrm{i}\mathrm{J} + \hat{\lambda}_{\mathsf{s}1}\hat{\lambda}_{\mathsf{s}2}\mathrm{K}^4 - \hat{\lambda}_{\mathsf{s}1}\mathrm{K}^2, \quad \textit{for } \mathsf{s} = 0, 1.$$

*All inverses, fractions and multiplications are here intended elementwise.*

## 8.1 Y bound

In the continuation case, we want to compute $Y$ such that

$$\max_{s\in[0,1]} \|A_s H_s(\hat{x}_s)\|_X \leqslant Y.$$

We denote $A_s H_s(\hat{x}_s) = (g_1(s), g_2(s), g_3(s)) \in X = \mathbb{R}^3 \times \ell^1_{\nu_2} \times \mathcal{L}^1_\nu$ for $s \in [0, 1]$.

In a similar way as in Lemma 8.2 it holds, for $i = 1, 2, 3$, that

$$\max_{s \in [0,1]} \|g_i(s)\|_{X_i} \leq \frac{1}{8} \left\| \left( \max_{s \in [0,1]} |g_k''(s)| \right)_{k \in I_i} \right\|_{X_i} \\ + \left\| (\max\{|(g_i(0))_k|, |(g_i(1))_k|\})_{k \in I_i} \right\|_{X_i}. \tag{70}$$

Applying the approach presented in Section 7.1 we can compute all nonvanishing components of $g_i(0)$ and $g_i(1)$, thus making the second term straightforward to bound explicitly. We are left with having to bound the second order derivative. It holds, by linearity in $s$ of $\hat{x}_s$ and $A_s$, that

$$g''(s) = D_s D_s (A_s H_s(\hat{x}_s)) = 2 A_\Delta DH(\hat{x}_s) x_\Delta + A_s D_s D_s H(\hat{x}_s)$$
$$= 2 A_\Delta DH(\hat{x}_s) x_\Delta + A_s DDH(\hat{x}_s) x_\Delta x_\Delta.$$

By construction, $x_\Delta \in \Pi^{(K)} X$, while $DH(\hat{x}_s)$ acts on $x_\Delta$ as a convolution. Therefore $DH(\hat{x}_s) x_\Delta \in \Pi^{(2K)} X$. By Lemma 4.6, thanks to $A_\Delta$ being eventually diagonal of center dimension $K$, we find that $A_\Delta DH(\hat{x}_s) x_\Delta \in \Pi^{(2K)} X$. In the same way we can prove that

$$A_s DDH(\hat{x}_s) x_\Delta x_\Delta \in \Pi^{(2K)} X.$$

Thus, we can bound $\|D_s D_s (A_s H_s(\hat{x}_s))\|_c$ by a computation where we replace $\hat{x}_s$ by the *interval*

$$x_s = \hat{x}_0 + [0, 1] \cdot x_\Delta, \tag{71}$$

interpreted elementwise, and all the computations involving $\hat{x}_s$ are performed using interval arithmetic. Thus, both terms in (70) can be bounded explicitly.

## 8.2  $Z_0$ bound

Applying the splitting (60), the $Z_0$ bound is such that

$$Z_0 \geq \max_s \|I - A_s A_s^\dagger\|_{B(X,X)}.$$

Using the notation (31) we apply Lemma 8.2 to obtain, for $i, j = 1, 2, 3$,

$$\|(I - A_s A_s^\dagger)_{ij}\|_{ij} \leq$$
$$\left\| \left( \max\{|((I - A_0 A_0^\dagger)_{ij})_{k_1 k_2}|, |((I - A_1 A_1^\dagger)_{ij})_{k_1 k_2}|\} \right)_{k_1 \in I_1, k_2 \in I_2} \right\|_{ij}$$
$$+ \frac{1}{4} \left\| \left( \sup_s |(A_\Delta A_\Delta^\dagger)_{ij})_{k_1 k_2}| \right)_{k_1 \in I_1, k_2 \in I_2} \right\|_{ij}.$$

The first term is the elementwise maximum over $I - A_0 A_0^\dagger$ and $I - A_1 A_1^\dagger$. We computed these terms and their norm separately in Section 7.3. Here we simply need to take the elementwise maximum of their absolute values *before* taking the norm, but it is otherwise a straightforward extension of the computation in Section 7.3. We are left to compute $\|(A_\Delta A_\Delta^\dagger)_{ij}\|_{ij}$, where the maximum over $s$ is dropped, because the argument does not depend on $s$. The product $A_\Delta A_\Delta^\dagger$ is the product of two eventually diagonal operator of center dimension $K$, thus it is again an eventually diagonal operator of center dimension $K$. It then holds that

$$\|I - A_\Delta A_\Delta^\dagger\| \leq \max\{\|\Pi^{(K)} - A_\Delta^F A_\Delta^{\dagger F}\|, \|\Pi_{(K)}^\infty - A_\Delta^I A_\Delta^{\dagger I}\|\}$$

The norm of the finite center can be computed directly, hence we will concentrate on the bound of the tail operator.

Let us recall, from the definitions of $A_i$ and $A_i^\dagger$, with $i = 0, 1$, that

$$A_\Delta^{\dagger I} = A_0^{\dagger I} - A_1^{\dagger I} = \begin{pmatrix} 0 & 0 & 0 \\ 0 & \hat{\lambda}_{2\Delta} K^4 & 0 \\ 0 & 0 & (\hat{\lambda}_1 \hat{\lambda}_2)_\Delta K^4 - \hat{\lambda}_{1\Delta} K^2 \end{pmatrix} \Pi_{(K)}^\infty, \qquad (72)$$

and $A_\Delta^I$ is defined in (69). It follows from Equations (72) and (69) that

$$A_\Delta^I A_\Delta^{\dagger I} = \begin{pmatrix} 0 & 0 & 0 \\ 0 & B_1 & 0 \\ 0 & 0 & B_2 \end{pmatrix} \Pi_{(K)}^\infty$$

with

$$B_1 = \frac{(\hat{\lambda}_{2\Delta} K^4)(\hat{\lambda}_{2\Delta} K^4)}{(\hat{\lambda}_{21} K^4 - K^2)(\hat{\lambda}_{20} K^4 - K^2)},$$
$$B_2 = \frac{((\hat{\lambda}_1 \hat{\lambda}_2)_\Delta K^4 - \hat{\lambda}_{1\Delta} K^2)((\hat{\lambda}_1 \hat{\lambda}_2)_\Delta K^4 - \hat{\lambda}_{1\Delta} K^2)}{(iJ + \hat{\lambda}_{11} \hat{\lambda}_{21} K^4 - \hat{\lambda}_{11} K^2)(iJ + \hat{\lambda}_{10} \hat{\lambda}_{20} K^4 - \hat{\lambda}_{10} K^2)}.$$

With this format, we could apply directly the results of Section 6, but we prefer instead to avoid computing bounds having diagonal terms of order 8. We therefore simplify the computations by considering

$$\|B_1\Pi^\infty_{(K)}\|_{22} \leq \left\|\frac{\hat{\lambda}_{2\Delta}K^4}{\hat{\lambda}_{21}K^4-K^2}\Pi^\infty_{(K)}\right\|_{22} \left\|\frac{\hat{\lambda}_{2\Delta}K^4}{\hat{\lambda}_{20}K^4-K^2}\Pi^\infty_{(K)}\right\|_{22},$$

$$\|B_2\Pi^\infty_{(K)}\|_{33} \leq \left\|\frac{(\hat{\lambda}_1\hat{\lambda}_2)_\Delta K^4-\hat{\lambda}_{1\Delta}K^2}{\mathrm{i}J+\hat{\lambda}_{11}\hat{\lambda}_{21}K^4-\hat{\lambda}_{11}K^2}\Pi^\infty_{(K)}\right\|_{33} \left\|\frac{(\hat{\lambda}_1\hat{\lambda}_2)_\Delta K^4-\hat{\lambda}_{1\Delta}K^2}{\mathrm{i}J+\hat{\lambda}_{10}\hat{\lambda}_{20}K^4-\hat{\lambda}_{10}K^2}\Pi^\infty_{(K)}\right\|_{33}.$$

In all these four cases, we can apply the approach presented in Section 6 to bound the norms of the rational tail operators.

This concludes the computation of the $Z_0$ bound in the continuation setting.

### 8.3 $Z_1$ bound

Applying Lemma 8.2, and using the notation

$$A_s(A_s^\dagger - DH_s(\hat{x}_s)) = \begin{pmatrix} B_{11}(s) & B_{12}(s) & B_{13}(s) \\ B_{21}(s) & B_{22}(s) & B_{23}(s) \\ B_{31}(s) & B_{32}(s) & B_{33}(s) \end{pmatrix}$$

we can write

$$\max_{s\in[0,1]} \sup_{c\in B_r(0)} \|B_{ij}(s)\|_{ij} \leq$$

$$\left\|(\max\{|(B_{ij}(1))_{k_1k_2}|, |(B_{ij}(0))_{k_1k_2}|\})_{k_1\in I_1, k_2\in I_2}\right\|_{ij}$$
$$+ \frac{1}{8}\left\|\max_{s\in[0,1]}|(B''_{ij}(s))_{k_1k_2}|\right\|_{ij} \leq Z_1 r,$$

for $i,j = 1,2,3$. The bound at the end points, $s = 0,1$, follows from the procedure presented in Section 7.4 taking the elementwise maxima before taking the norms, as in Section 8.2. Therefore it remains to bound the components of

$$B''(s) = D_s D_s \left(A_s(A_s^\dagger - DH_s(\hat{x}_s))\right).$$

By linearity of $A_s$ and $A_s^\dagger$ in $s$, it holds that

$$D_s D_s \left(A_s(A_s^\dagger - DH_s(\hat{x}_s))\right) \tag{73}$$
$$= 2(A_\Delta D_s\left(A_s^\dagger - DH_s(\hat{x}_s)\right) + A_s D_s D_s \left(A_s^\dagger - DH_s(\hat{x}_s)\right) \tag{74}$$
$$= 2A_\Delta A_\Delta^\dagger - 2A_\Delta D^2 H_s(\hat{x}_s)x_\Delta - A_s D^3 H_s(\hat{x}_s)x_\Delta x_\Delta.$$

We already discussed a bound for $A_\Delta A_\Delta^\dagger$ in Section 8.2. For the second term $A_\Delta D^2 H_s(\hat{x}_s) x_\Delta$, we recall the direct computation of $D^2 H_s$ in (65). Furthermore, we have $\hat{x}_s, x_\Delta \in \Pi^{(K)} X$. Since $A_\Delta$ is an eventually diagonal operator, it holds that

$$A_\Delta D^2 H_s(\hat{x}_s) x_\Delta = A_\Delta^F \Pi^{(K)} D^2 H_s(\hat{x}_s) x_\Delta + A_\Delta^I \Pi_{(K)}^\infty D^2 H_s(\hat{x}_s) x_\Delta.$$

In the first term, we now replace, as in Section 8.1, $\hat{x}_s$ by the interval $x_s$ defined in (71). Since we know that $D^2 H_s(\hat{x}_s)$ has a finite number of nonvanishing (interval) components, we can explicitly compute all the nonvanishing components of the first term and compute a bound for $\|A_\Delta^F \Pi^{(K)} D^2 H_s(\hat{x}_s) x_\Delta\|_c$.

For the second term, it holds that

$$A_\Delta^I \Pi_{(K)}^\infty D^2 H_s(\hat{x}_s) x_\Delta = \begin{pmatrix} 0 & 0 & 0 \\ 0 & A_{\Delta 1}^I d_{22}(\hat{x}_s, x_\Delta) & 0 \\ 0 & A_{\Delta 2}^I d_{32}(\hat{x}_s, x_\Delta) & A_{\Delta 2}^I d_{33}(\hat{x}_s, x_\Delta) \end{pmatrix}.$$

With the notation introduced in (69), we can write

$$A_{\Delta 1}^I d_{22}(\hat{x}_s, x_\Delta) = $$
$$\frac{\hat{\lambda}_{\Delta 2} K^4}{(\hat{\lambda}_{20} K^4 - K^2)(\hat{\lambda}_{21} K^4 - K^2)} \Pi_{(K)}^\infty d_{22}(\hat{x}_s, x_\Delta),$$

$$A_{\Delta 2}^I d_{32}(\hat{x}_s, x_\Delta) = $$
$$\frac{(\hat{\lambda}_1 \hat{\lambda}_2)_\Delta K^4 - \hat{\lambda}_{1\Delta} K^2}{(iJ + \hat{\lambda}_{10} \hat{\lambda}_{20} K^4 - \hat{\lambda}_{10} K^2)(iJ + \hat{\lambda}_{11} \hat{\lambda}_{21} K^4 - \hat{\lambda}_{11} K^2)} \Pi_{(K)}^\infty d_{32}(\hat{x}_s, x_\Delta),$$

$$A_{\Delta 2}^I d_{33}(\hat{x}_s, x_\Delta) = $$
$$\frac{(\hat{\lambda}_1 \hat{\lambda}_2)_\Delta K^4 - \hat{\lambda}_{1\Delta} K^2}{(iJ + \hat{\lambda}_{10} \hat{\lambda}_{20} K^4 - \hat{\lambda}_{10} K^2)(iJ + \hat{\lambda}_{11} \hat{\lambda}_{12} K^4 - \hat{\lambda}_{11} K^2)} \Pi_{(K)}^\infty d_{33}(\hat{x}_s, x_\Delta).$$

We apply a splitting of the fraction similar to the one used in Section 8.2 to simplify the computations, thereby retrieving

$$\|A^I_{\Delta 1} d_{22}(\hat{x}_s, x_\Delta)\|_{22} =$$
$$\left\|\Pi^\infty_{(K)} \frac{d_{22}(\hat{x}_s, x_\Delta)}{\hat{\lambda}_{20}K^4 - K^2}\right\|_{22} \left\|\Pi^\infty_{(K)} \frac{\hat{\lambda}_{2\Delta}K^4}{\hat{\lambda}_{21}K^4 - K^2}\right\|_{22},$$

$$\|A^I_{\Delta 2} d_{32}(\hat{x}_s, x_\Delta)\|_{32} =$$
$$\left\|\Pi^\infty_{(K)} \frac{d_{32}(\hat{x}_s, x_\Delta)}{iJ + \hat{\lambda}_{10}\hat{\lambda}_{20}K^4 - \hat{\lambda}_{10}K^2}\right\|_{32} \left\|\Pi^\infty_{(K)} \frac{(\hat{\lambda}_1\hat{\lambda}_2)_\Delta K^4 - \hat{\lambda}_{1\Delta}K^2}{iJ + \hat{\lambda}_{11}\hat{\lambda}_{21}K^4 - \hat{\lambda}_{11}K^2}\right\|_{32},$$

$$\|A^I_{\Delta 2} d_{33}(\hat{x}_s, x_\Delta)\|_{33} =$$
$$\left\|\Pi^\infty_{(K)} \frac{d_{33}(\hat{x}_s, x_\Delta)}{iJ + \hat{\lambda}_{10}\hat{\lambda}_{20}K^4 - \hat{\lambda}_{10}K^2}\right\|_{33} \left\|\Pi^\infty_{(K)} \frac{(\hat{\lambda}_1\hat{\lambda}_2)_\Delta K^4 - \hat{\lambda}_{1\Delta}K^2}{iJ + \hat{\lambda}_{11}\hat{\lambda}_{21}K^4 - \hat{\lambda}_{11}K^2}\right\|_{33}.$$

The second norm in each bound can be computed by applying the procedure presented in Section 6, but the first norms are not yet in that form, because $d_{ij}(\hat{x}_s, x_\Delta)$ are not diagonal operators. Using the expressions for $d_{ij}$ introduced in Section 7.5, we estimate

$$\left\|\Pi^\infty_{(K)} \frac{d_{22}(\hat{x}_s, x_\Delta)}{\hat{\lambda}_{20}K^4 - K^2}\right\|_{22} \leqslant$$
$$\left\|\Pi^\infty_{(K)} \frac{K^4}{\hat{\lambda}_{20}K^4 - K^2}\right\|_{22} |\lambda_{2\Delta}| + \left\|\Pi^\infty_{(K)} \frac{2K}{\hat{\lambda}_{20}K^4 - K^2}\right\|_{22} \|y_\Delta\|_{\ell^1_{\nu_2}},$$

and

$$\left\|\Pi^\infty_{(K)} \frac{d_{32}(\hat{x}_s, x_\Delta)}{iJ + \hat{\lambda}_{10}\hat{\lambda}_{20}K^4 - \hat{\lambda}_{10}K^2}\right\|_{32} \leqslant$$
$$\left\|\Pi^\infty_{(K)} \frac{2K}{iJ + \hat{\lambda}_{10}\hat{\lambda}_{20}K^4 - \hat{\lambda}_{10}K^2}\right\|_{33} \|\lambda_{1s} z_\Delta + \lambda_{1\Delta} z_s\|_{\mathcal{L}^1_\nu},$$

and

$$\left\|\Pi_{(K)}^\infty \frac{d_{33}(\hat{x}_s, x_\Delta)}{\mathrm{i}J + \hat{\lambda}_{10}\hat{\lambda}_{20}K^4 - \hat{\lambda}_{10}K^2}\right\|_{33} \leqslant$$

$$\left\|\Pi_{(K)}^\infty \frac{K^4}{\mathrm{i}J + \hat{\lambda}_{10}\hat{\lambda}_{20}K^4 - \hat{\lambda}_{10}K^2}\right\|_{33} |\lambda_{2s}\lambda_{1\Delta} + \lambda_{1s}\lambda_{2\Delta}|$$

$$+ \left\|\Pi_{(K)}^\infty \frac{K^2}{\mathrm{i}J + \hat{\lambda}_{10}\hat{\lambda}_{20}K^4 - \hat{\lambda}_{10}K^2}\right\|_{33} |\lambda_{1\Delta}|$$

$$+ \left\|\Pi_{(K)}^\infty \frac{2K}{\mathrm{i}J + \hat{\lambda}_{10}\hat{\lambda}_{20}K^4 - \hat{\lambda}_{10}K^2}\right\|_{33}$$

$$+ \|(E_\ell y_s + a_s z_s)\lambda_{1\Delta} + \lambda_{1s}(a_\Delta z_s + a_s z_\Delta + E_\ell y_\Delta)\|_{\mathcal{L}_\nu^1}.$$

All these terms can either be bounded by direct interval arithmetic computation (since both $x_\Delta$ and the interval-valued $x_s \ni \hat{x}_s$ have finitely many nonvanishing components only) or via the bounds on rational tail operators in Section 6.

We will now bound the final term in (73), $A_s D^3 H_s(\hat{x}_s) x_\Delta x_\Delta$. This term first requires the computation of $D^3 H_s(\hat{x}_s) x_\Delta x_\Delta$, that is

$$D^3 H(\hat{x}_s) x_\Delta x_\Delta = \begin{pmatrix} 0 & 0 & 0 \\ 0 & 0 & 0 \\ d_{31}^3(\hat{x}_s, x_\Delta) & d_{32}^3(\hat{x}_s, x_\Delta) & d_{33}^3(\hat{x}_s, x_\Delta) \end{pmatrix},$$

with

$$d_{31}^3(\hat{x}_s, x_\Delta) = \left[d_{311}^3(\hat{x}_s, x_\Delta), 2K^4 z_\Delta \lambda_{1\Delta}, d_{313}^3(\hat{x}_s, x_\Delta)\right],$$
$$d_{311}^3(\hat{x}_s, x_\Delta) = 2K^4 z_\Delta \lambda_{2\Delta} - \mathrm{i}K(4E_\ell y_\Delta * z_\Delta + 4a_\Delta z_\Delta * z_s + 2a_s z_\Delta * z_\Delta)$$
$$d_{313}^3(\hat{x}_s, x_\Delta) = -\mathrm{i}K(4\lambda_{1\Delta} z_\Delta * z_s + 2\lambda_{1s} z_\Delta * z_\Delta),$$
$$d_{32}^3(\hat{x}_s, x_\Delta) = -4\mathrm{i}K\lambda_{1\Delta} C^{\mathcal{L}_\nu^1}(z_\Delta) E_\ell,$$
$$d_{33}^3(\hat{x}_s, x_\Delta) = 2K^4 \lambda_{1\Delta} \lambda_{2\Delta}$$
$$\qquad - 4\mathrm{i}KC^{\mathcal{L}_\nu^1}(\lambda_{1s} a_\Delta z_\Delta + \lambda_{1\Delta}(E_\ell y_\Delta + a_\Delta z_s + a_s z_\Delta)).$$

The norm of $A_s D^3 H(\hat{x}_s) x_\Delta x_\Delta$ can be bounded by considering $A_s$ as the interval $[A_0, A_1]$, interpreted elementwise and using the triangular inequality:

$$\|A_s D^3 H(\hat{x}_s) x_\Delta x_\Delta\| \leqslant \|A_s^I D^3 H(\hat{x}_s) x_\Delta x_\Delta\| + \|A_s^F D^3 H(\hat{x}_s) x_\Delta x_\Delta\|.$$

Here we can bound the finite part with Corollary 4.13, since all nonvanishing terms can be enclosed using interval arithmetic. The infinite tail gives

$$A_s^I D^3 H(\hat{x}_s) x_\Delta x_\Delta = \begin{pmatrix} 0 & 0 & 0 \\ 0 & 0 & 0 \\ 0 & A_{s2}^I d_{32}^3 & A_{s1}^I d_{33}^3 \end{pmatrix},$$

and

$$\|A_{s2}^I d_{32}^3\|_{32} \leqslant 4\|A_{s2}^I \mathrm{K}\|_{32} \|\lambda_{1\Delta} z_\Delta\|_{\mathcal{L}_\nu^1},$$
$$\|A_{s2}^I d_{33}^3\|_{33} \leqslant 2\|A_{s2}^I \mathrm{K}^4\|_{33} |\lambda_{1\Delta}\lambda_{2\Delta}| + 4\|A_{s2}^I \mathrm{K}\|_{33} \|4$$
$$+ \lambda_{1s} a_\Delta z_\Delta + \lambda_{1\Delta}(E_\ell y_\Delta + a_\Delta z_s + a_s z_\Delta)\|_{\mathcal{L}_\nu^1}.$$

These terms can be bounded either by direct computation or by the technique presented in Section 6.

## 8.4 $Z_2$ bound

For the $Z_2$ bound in the continuation case, we proceed in the same way as in Section 7.5. By applying the mean value theorem we have

$$\max_{s \in [0,1]} \sup_{b,c \in B_r(0)} \|A_s \left(DH_s(\hat{x}_s) - DH_s(\hat{x}_s + rb)rc\right)\| \leqslant$$
$$\max_{s \in [0,1]} \sup_{b,c \in B_1(0)} \|A_s D^2 H_s(\hat{x}_s + Rb)c\| r^2.$$

Following the strategy presented in Section 7.5 and using the convexity of the norm, we bound

$$\max_{s \in [0,1]} \sup_{\bar{b},\bar{c} \in B_1(0)} \sup_{z \in B_1(0)} \left\|\left(A_s D^2 H_s(\hat{x}_s + zr)\bar{b}\bar{c}\right)_n\right\| r^2 \leqslant$$
$$\max_{A=\{A_0,A_1\}} \sup_{\bar{b},\bar{c} \in B_1(0)} \sup_{z \in B_1(0)} \left\|\left(A D^2 H_s(\hat{x}_s + zr)\bar{b}\bar{c}\right)_n\right\| r^2,$$

and we can apply the same process as in Section 7.5 of estimating "monomial by monomial", with the additional step that

$$\|\hat{x}_s\|_c \leqslant \max\{\|\hat{x}_0\|_c, \|\hat{x}_1\|_c\},$$

where the component norm $\|\cdot\|_c$ is defined in (20).

# 9 Conclusion

We have until now discussed how to validate already computed numerical solutions. Before presenting the results, we explain how the first numerical solution $\hat{x} = (\hat{\lambda}_1, \hat{\lambda}_2, \hat{a}, \hat{y}, \hat{z})$ is computed. We consider the stationary problem (13) projected onto a finite dimensional subset $\Pi^{(K)} \ell^1_{\nu_2} \subset \ell^1_{\nu_2}$. In this framework, Equation (13) becomes a finite dimensional system which we can easily solve numerically for $y$ at any given $\lambda_2 = \gamma$. Then, it is straightforward to compute all the eigenvalues $(\mu_1, \ldots \mu_{2K+1})$ associated to (13) at a given $\gamma$ and $y$. This procedure does not need to be rigorous. We look for a value $\gamma = \lambda_2$ and a corresponding numerical approximate solution $\hat{y} \in \tau^{(K)} \ell^1_{\nu_2}$ of Equation (13) such that a pair of eigenvalues $\mu$ and $\bar{\mu}$ crosses the imaginary axes. Let $v$ be the eigenvector associated to $\mu$. We set

$$\hat{z}_{jk} = \begin{cases} v_k & \text{if } j = 1, |k| \leqslant K, \\ \bar{v}_{-k} & \text{if } j = -1, |k| \leqslant K, \\ 0 & \text{otherwise.} \end{cases}$$

We also set $\hat{a} = 0$ as a first approximation, because we assume to be very close to a Hopf bifurcation, and we set $\hat{\lambda}_1$ to $|\text{Im}(\mu)|^{-1}$. This construction forms a first numerical approximation for $\hat{x}$. Applying Newton's method allows us to increase the accuracy of our initial approximation. In Figure 1 both the stationary solution $\hat{y}$ at $\gamma \approx \gamma^*$ close to the Hopf bifurcation and the rescaled time dependent perturbation $\hat{z}$ are depicted.

Once this approximation has been constructed, the predictor-corrector algorithm, as presented for example in [3] and [10], allows us to numerically follow the solution branch. Then, with the method presented in this article, we rigorously validate each segment and follow the solution branch. The validated segment where the sign of $a$ is proven to change is where the Hopf bifurcation happens, based on the arguments in Remarks 2.6 and 2.7.

The MATLAB script `script_HopfKS.m`, which can be found at [11], performs the necessary computation to conclude that we have validated the existence of a Hopf bifurcation in the Kuramoto-Sivashinsky equation from the stationary solution depicted in Figure 1 at parameter value $\lambda_2^* = \gamma^* = 0.298776358114 + [-r, r]$ with the validation radius $r = 8.63 \cdot 10^{-11}$. For the validation we used 20 temporal modes and 30 spatial modes: $K = (20, 30)$. The weights used for the $\mathcal{L}^1_\nu$ and $\ell^1_{\nu_2}$ space are $\nu = (1.1, 1.1)$. The segment has been followed from $a = -10^{-8} < 0$ to $a = 5.3898 \cdot 10^{-4} > 0$. With a step size of $10^{-6}$, this validation required 539 steps and reached $\lambda_2 = 0.29877636757 + [-r, r]$, with $r = 1.6 \cdot 10^{-10}$, which is at least $8 \cdot 10^{-9}$ above the bifurcation value $\lambda_2^*$.